\documentclass[12pt,sort&compress]{elsarticle}

\usepackage{bm}
\usepackage{comment}

\usepackage{mathrsfs}
\usepackage{amsmath}
\usepackage{amsfonts}
\usepackage{amsthm}
\usepackage{color}
\usepackage{xcolor}

\usepackage{caption}
\usepackage{subcaption}

\usepackage{float}
\usepackage{url}
\usepackage{multicol}
\usepackage[top=1in,bottom=1in,left=1in,right=1in]{geometry}
\usepackage[pdfborder={0 0 0},colorlinks,allcolors=blue]{hyperref}
\usepackage{txfonts}
\usepackage{setspace}
\usepackage{enumitem}
\usepackage{lipsum} 
\usepackage{siunitx,booktabs}
\usepackage{nth}
\usepackage{accents} 

\usepackage{makecell}
\usepackage{appendix}

\theoremstyle{definition} 

\allowdisplaybreaks

\newcommand{\cmmnt}[1]{\ignorespaces}

\usepackage{listings}
\definecolor{light-gray}{gray}{0.95}
\begin{document}
\begin{frontmatter}

\title{Shape optimization of non-matching isogeometric shells with moving intersections}

\author[ucsdmae]{Han Zhao} 
\author[ucsdmae]{John T. Hwang} 
\author[ucsdse,ucsdmae]{J. S. Chen\corref{cor1}} 
\ead{jsc137@ucsd.edu}
\cortext[cor1]{Corresponding author}

\affiliation[ucsdmae]{organization={Department of Mechanical and Aerospace Engineering, University of California San Diego},
            addressline={9500 Gilman Drive}, 
            city={La Jolla},
            postcode={92122}, 
            state={CA},
            country={USA}}

\affiliation[ucsdse]{organization={Department of Structural Engineering, University of California San Diego},
            addressline={9500 Gilman Drive}, 
            city={La Jolla},
            postcode={92122}, 
            state={CA},
            country={USA}}


\begin{abstract}
    While shape optimization using isogeometric shells exhibits appealing features by integrating design geometries and analysis models, challenges arise when addressing computer-aided design (CAD) geometries comprised of multiple non-uniform rational B-splines (NURBS) patches, which are common in practice. The intractability stems from surface intersections within these CAD models. In this paper, we develop an approach for shape optimization of non-matching isogeometric shells incorporating intersection movement. Separately parametrized NURBS surfaces are modeled using Kirchhoff--Love shell theory and coupled using a penalty-based formulation. The optimization scheme allows shell patches to move without preserving relative location with other members during the shape optimization. This flexibility is achieved through an implicit state function, and analytical sensitivities are derived for the relative movement of shell patches. The introduction of differentiable intersections expands the design space and overcomes challenges associated with large mesh distortion, particularly when optimal shapes involve significant movement of patch intersections in physical space. Throughout optimization iterations, all members within the shell structures maintain the NURBS geometry representation, enabling efficient integration of analysis and design models. The optimization approach leverages the multilevel design concept by selecting a refined model for accurate analysis from a coarse design model while maintaining the same geometry. We adopt several example problems to verify the effectiveness of the proposed scheme and demonstrate its applicability to the optimization of the internal stiffeners of an aircraft wing. 
\end{abstract}


\end{frontmatter}

\section{Introduction}\label{sec:intro}

The Kirchhoff--Love shell model requires $C^1$ continuous basis functions, isogeometric analysis (IGA) \cite{Hughes05a} using NURBS \cite{piegl2012nurbs} basis functions is perfectly suited for the solution of Kirchhoff–Love shells. Moreover, IGA offers a unified framework for seamless integration of CAD models and analysis models, circumventing the cumbersome process of finite element (FE) mesh generation \cite{Hardwick2005}. Comprehensive research on isogeometric Kirchhoff--Love shell model is conducted in \cite{Kiendl2009, Kiendl2011, Kiendl2015, NguyenThanh2011, Casquero2017, Casquero2020, magisano2024large}. Various applications including heart valve leaflets \cite{hsu2014fluid, Kamensky2015, Kamensky2017b, Kamensky2021, neighbor2023leveraging}, wind turbines \cite{Herrema2019a, herrema2019penalty, johnson2020isogeometric}, composite materials \cite{schulte2020isogeometric, mohammadi2022isogeometric}, and aerospace structures \cite{hirschler2019isogeometric, zhao2022open, hao2023isogeometric} have demonstrated the capability of the Kirchhoff--Love shell theory with isogeometric discretization. 

Multiple patches are typically required to model complex, realistic shell structures using NURBS surfaces. To make the CAD geometries with multiple patches directly available for structural analysis, coupling between adjacent NURBS patches to maintain the displacement and rotational continuities across patch intersections becomes essential. The bending strip method \cite{Bazilevs10c} and the kinematic constraints  \cite{Duong2017} have been proposed for coupling NURBS surfaces with conforming discretizations. For Kirchhoff--Love shells with non-matching intersections, a series of coupling techniques have been explored, including mortar methods \cite{BRIVADIS2015292, horger2019hybrid, hirschler2019dual}, Nitsche-type methods \cite{guo2015nitsche, Guo2018, guo2021isogeometric, Benzaken2021, wang2022isogeometric, yu2023isogeometric}, projected super-penalty methods \cite{Coradello2021, Coradello2021a}, penalty methods \cite{herrema2019penalty,leonetti2020robust,proserpio2022penalty,zhao2022open, guarino2024interior}, and embedded surfaces methods \cite{hirschler2019embedded, hao2023isogeometric}.

A well-designed shell structure features excellent performance by distributing load through membrane forces while minimizing bending moments \cite{bletzinger1993form}, with the mechanical characteristics significantly affected by its shape. Consequently, shape optimization plays a critical role in the development of novel shell structures. The unified model between geometric design and structural analysis in IGA renders particular advantages for shell shape optimization \cite{wall2008isogeometric, cho2009isogeometric, ha2010numerical, qian2010full, li2011isogeometric, azegami2013shape}, with many superior designs such as composite shells \cite{nagy2013isogeometric, hao2018integrated}, wind turbine blades \cite{herrema2017framework, Herrema2019a}, and stiffened thin-wall structures \cite{hirschler2019isogeometric, hao2023isogeometric}. Traditionally, shape optimization relies on the finite element method (FEM) with parametric models \cite{ding1986shape}. However, the classical FEM-based approach encounters difficulties in the precise representation of the updated geometry and accurate solution of structural behavior \cite{imam1982three}. Shape optimization using IGA addresses these challenges by directly performing structural analysis on the design model, circumventing the tedious intermediate steps as in the FEM-based approach. The geometric error is eliminated and the continuity of the geometry is preserved by adjusting the coordinates of the control points during the optimization process. Nonetheless, updating complex geometries with multiple NURBS patches necessitates additional efforts to represent surface intersections accurately. \cite{hirschler2019embedded, hao2023isogeometric} used embedded surfaces in an extruded free-form deformation (FFD) block \cite{Sederberg1986} to impose shape modifications. \cite{zhao2024automated} employed the FFD idea in conjunction with Lagrange extraction \cite{Schillinger2016} to perform shape optimization for the non-matching shell patches while maintaining the intersection geometries. However, challenges remain with these methods. The ``master'' surfaces need to be identified, and extrusion to a 3D B-spline block is required in the embedded surfaces method, where the latter approach may lead to substantial distortion of elements when surface intersections undergo large movement.

In this paper, we propose a shape optimization method for isogeometric Kirchhoff--Love shell structures consisting of a stack of NURBS surfaces with moving intersections. The control points of all shell patches are optimized directly without additional effort, while the locations of the surface intersections are updated accordingly through an implicit relation between the control points and intersections' parametric coordinates. Hence, relative movement between shell patches is made without distorting the shell elements. The framework in \cite{zhao2022open} is employed to couple the non-matching shell patches using a penalty-based formulation. The coupling method involves creating a topologically 1D quadrature mesh in the parametric space to integrate the penalty energy at the intersection positions. Throughout the optimization process, parametric coordinates of intersections are solved accordingly when updating the geometry of shell patches. Sensitivities of the implicit relation and penalty residual with respect to intersections' parametric coordinates are derived to obtain the total derivative of the optimization problem. In this approach, separately modeled NURBS patches can move smoothly relative to other intersecting patches as long as the intersections exist. The quality of the shell element is insensitive to the large movement of the intersection, avoiding ill-conditioning in the discrete system. Additionally, all NURBS surfaces can be parametrized without distinguishing between ``master'' and ``slave'' surfaces and performing 3D solid extrusion. This approach is particularly beneficial for the design of internal structures of an aircraft wing, where the placement of the internal structures can be determined straightforwardly without compromising the element quality of the outer surfaces. We demonstrate the effectiveness of this approach through innovative designs for the internal structures of an electric vertical takeoﬀ and landing (eVTOL) aircraft wing.

The remainder of the paper is organized as follows. Section \ref{sec:coupling-KL-shell} provides a review of the non-matching isogeometric Kirchhoff--Love shells coupling formulations and algorithms. Section \ref{sec:shape-opt-nonmatching-shells} presents the shape optimization approach for non-matching shells with moving intersections, and derives the total derivative with respect to the design variables. Section \ref{sec:implementation} discusses the implementation details and the associated numerical procedures for the optimization framework. Two benchmark problems with reference solutions are used to validate the shape optimization approach in Section \ref{sec:benchmark}, followed by a demonstration of its applicability to the internal structure shape optimization of an eVTOL wing in Section \ref{sec:applications}. Finally, Section \ref{sec:conclusion} draws the concluding remarks of the proposed optimization approach.


\section{Non-matching coupling of Kirchhoff--Love shell} \label{sec:coupling-KL-shell}
Structural analysis is crucial for the evaluation of the structural performance and sensitivity calculation for shape optimization. In this work, shell structures are modeled using the Kirchhoff--Love shell theory discretized by NURBS basis functions with higher order continuity. Under this framework, separate shell patches in the CAD geometry are coupled using a penalty-based formulation.

\subsection{Basic Kirchhoff--Love shell formulation}\label{subsec:basic-KL-shell}
This section only provides an overview to lay the foundation for the subsequent optimization approach. In the Kirchhoff--Love shell theory \cite{Kiendl2011}, the 3D shell continuum is represented by its mid-surface, which can be parametrized in a 2D space using coordinates $\bm{\xi}=\{\xi_1, \xi_2\}$. We denote the geometry of the mid-surface in the reference configuration as $\mathbf{X}(\bm{\xi})$ and the deformed configuration as $\mathbf{x}(\bm{\xi})$. The displacement field of the mid-surface is given by 
\begin{align}
    \mathbf{x}(\bm{\xi}) = \mathbf{X}(\bm{\xi}) + \mathbf{u}(\bm{\xi})\label{eq:kl-shell-deformed-config} \text{ .}
\end{align}
Covariant basis vectors of the mid-surface are defined as
\begin{align}
    \mathbf{A}_{\alpha} = \mathbf{X},_{\xi_\alpha} \qquad \text{and} \qquad \mathbf{a}_{\alpha} = \mathbf{x},_{\xi_\alpha} \label{eq:kl-shell-covariant-basis} \text{ ,}
\end{align}
where $(\cdot),_{\xi_\alpha} = \frac{\partial (\cdot)}{\partial \xi_\alpha}$ and $\alpha=\{1,2\}$.
Unit vectors that are normal to the mid-surface are given by
\begin{align}
    \mathbf{A}_{3} = \frac{\mathbf{A}_{1} \times \mathbf{A}_{2}}{\Vert \mathbf{A}_{1} \times \mathbf{A}_{2} \Vert} \qquad \text{and} \qquad \mathbf{a}_{3} = \frac{\mathbf{a}_{1} \times \mathbf{a}_{2}}{\Vert \mathbf{a}_{1} \times \mathbf{a}_{2} \Vert} \label{eq:kl-shell-normal-vectors} \text{ ,}
\end{align}
where $\Vert \cdot \Vert$ is the $L_2$ norm. With surface basis vectors in \eqref{eq:kl-shell-covariant-basis}, metric coefficients in both configurations are defined as
\begin{align}
    A_{\alpha \beta} = \mathbf{A}_{\alpha} \cdot \mathbf{A}_{\beta} \qquad \text{and} \qquad a_{\alpha \beta} = \mathbf{a}_{\alpha} \cdot \mathbf{a}_{\beta} \label{eq:kl-shell-metric-coeff} \text{ ,}
\end{align}
for $\alpha, \beta=\{1,2\}$, and curvature coefficients read as
\begin{align}
    B_{\alpha \beta} = \mathbf{A}_{\alpha,\,\xi_\beta} \cdot \mathbf{A}_3 = - \mathbf{A}_{\alpha} \cdot \mathbf{A}_{3,\,\xi_\beta} \qquad \text{and} \qquad b_{\alpha \beta} = \mathbf{a}_{\alpha,\,\xi_\beta} \cdot \mathbf{a}_3 = - \mathbf{a}_{\alpha} \cdot \mathbf{a}_{3,\,\xi_\beta} \label{eq:kl-shell-curvature-coeff} \text{ .}
\end{align}
The membrane strain tensor and curvature change tensor coefficients are formulated as
\begin{align}
    \varepsilon_{\alpha \beta} = \frac{1}{2}(a_{\alpha \beta} - A_{\alpha \beta}) \qquad \text{and} \qquad \kappa_{\alpha \beta} = B_{\alpha \beta} - b_{\alpha \beta} \label{eq:kl-shell-strain-coeff} \text{ .}
\end{align}
We employ the St. Venant--Kirchhoff material model in this paper with material tensor $\mathbf{C}$ to express normal forces and bending moments
\begin{align}
    \mathbf{n} = t\, \mathbf{C} : \bm{\varepsilon} \qquad \text{and} \qquad \mathbf{m} = \frac{t^3}{12}\, \mathbf{C} : \mathbf{\kappa} \label{eq:kl-shell-normal-force-bending-moments} \text{ .}
\end{align}
Using membrane strains and changes in curvature defined in \eqref{eq:kl-shell-strain-coeff} and associated force resultants in \eqref{eq:kl-shell-normal-force-bending-moments}, the virtual work of the Kirchhoff--Love shell read as
\begin{align}
   \delta W_{\text{s}} = \delta W^{\text{int}}_{\text{s}} - \delta W^{\text{ext}}_{\text{s}} = \int_S \delta \bm{\varepsilon} : \mathbf{n} + \delta \bm{\kappa} : \mathbf{m} \, \mathrm{d}S - \int_S \delta \mathbf{u} \cdot \mathbf{f} \, \mathrm{d}S \label{eq:kl-shell-virtual-work} \text{ ,}
\end{align}
where $S$ is the shell mid-surface and $\mathbf{f}$ is the external force acting on $S$, and $\delta W^{\text{int}}_{\text{s}}$ and $\delta W^{\text{ext}}_{\text{s}}$ represent the internal and external virtual work, respectively. \footnote{We use subscript ``$\text{s}$'' in symbols such as $\delta W^{\text{int}}_{\text{s}}$ and $\delta W^{\text{ext}}_{\text{s}}$ to denote the quantities on shell patches.}. 
A detailed derivation is presented in \cite[Section 3]{Kiendl2011}.

It is noted that the curvature coefficients in \eqref{eq:kl-shell-curvature-coeff} involve second-order derivatives of the displacements $\mathbf{u}$ and mid-surface geometry $\mathbf{X}$, basis functions with as least $C^1$ continuity on element boundaries is required. Discretization using NURBS basis functions automatically meets this requirement without additional treatment.

\subsection{Penalty coupling of shell patches}\label{subsec:shell-coupling}
Many complex shell structures comprise more than one NURBS patch. A coupling approach is needed for a collection of isogeometrically discretized shell patches to make them directly available for analysis. A penalty-based coupling formulation proposed by Herrema et al. \cite{herrema2019penalty} is employed in our current framework. The penalty energy preserves both displacement and rotational continuities on the intersection $\mathcal{L}$ between shell patch $S^\text{A}$ and $S^\text{B}$. The virtual work of the penalty energy is given by
\begin{equation}
\begin{aligned}
    \delta W^{\text{AB}}_{\text{pen}} &= \int_\mathcal{L} \alpha_\text{d}\, (\mathbf{u}^\text{A}-\mathbf{u}^\text{B})\cdot (\delta\mathbf{u}^\text{A}-\delta\mathbf{u}^\text{B}) \, d\mathcal{L}\\
    &+ \int_\mathcal{L} \alpha_\text{r} \, \left( (\mathbf{a}^\text{A}_3\cdot\mathbf{a}^\text{B}_3 - \mathbf{A}^\text{A}_3\cdot\mathbf{A}^\text{B}_3) (\delta \mathbf{a}^\text{A}_3\cdot\delta\mathbf{a}^\text{B}_3 - \delta\mathbf{A}^\text{A}_3\cdot\delta\mathbf{A}^\text{B}_3) \right. \\
    &\qquad \quad \left. +\,(\mathbf{a}^\text{A}_n\cdot\mathbf{a}^\text{B}_3 - \mathbf{A}^\text{A}_n\cdot\mathbf{A}^\text{B}_3) (\delta \mathbf{a}^\text{A}_n\cdot\delta\mathbf{a}^\text{B}_3 - \delta\mathbf{A}^\text{A}_n\cdot\delta\mathbf{A}^\text{B}_3) \right) \, \mathrm{d} \mathcal{L} \label{eq:penatly-virtual-work}\text{ ,}
\end{aligned}
\end{equation}
where $\mathbf{a}_3$ and $\mathbf{a}_n$ are normal and conormal vectors on the deformed configuration, while their counterparts in the reference configuration are denoted with uppercase letters. Computation of $\mathbf{a}_n$ is discussed in detail in Section \ref{subsec:isogeometric-discretization}. The scalar values $\alpha_d$ and $\alpha_r$ are penalty parameters for displacement and rotational continuities. These two parameters are constructed to account for material and geometric properties and are scaled by a problem-independent and dimensionless penalty coefficient $\alpha$
\begin{align}
    \alpha_{\text{d}} = \alpha \frac{Et}{h(1-\nu^2)} \qquad \text{ and } \qquad \alpha_{\text{r}} = \alpha \frac{Et^3}{12h(1-\nu^3)} \label{penalty-parameters} \text{ ,}
\end{align}
where $E$, $\nu$, and $t$ are Young's modulus, Poisson's ratio, and shell thickness, respectively. $h$ are averaged element length of shell patches $S^\text{A}$ and $S^\text{B}$. Details of the penalty formulation and coupling for composite shell structures can be found in \cite[Section 2]{herrema2019penalty}, where a wide range of effective penalty coefficients $\alpha$ was proposed. In this paper, we use the recommended value $\alpha=1000$ for all numerical examples.

With the virtual work of the shell patch in \eqref{eq:kl-shell-virtual-work} and penalty energy in \eqref{eq:penatly-virtual-work}, the total virtual work of two coupled shell patches $S^\text{A}$ and $S^\text{B}$ in the equilibrium state is expressed as
\begin{align}
    \delta W_{} = \delta W^{\text{A}}_{\text{s}} + \delta W^{\text{B}}_{\text{s}} + \delta W^{\text{AB}}_{\text{pen}} = 0 \text{ .} \label{eq:coupling-virtual-work}
\end{align}

\subsection{Shell coupling with isogeometric discretization} \label{subsec:isogeometric-discretization}
With the NURBS basis functions, Kirchhoff--Love shell geometry and displacement field are discretized isogeometrically. The position vector on the mid-surface of the shell patch in the reference configuration and the associated displacement vector are formulated as
\begin{align}
    \mathbf{X}(\bm{\xi}) = \sum\limits_{i=1}^n \hat{R}_{ip} (\bm{\xi}) \mathbf{P}_i = \hat{\mathbf{R}}(\bm{\xi}) \mathbf{P} \qquad \text{and} \qquad \mathbf{u}(\bm{\xi}) = \sum\limits_{i=1}^n \hat{R}_{ip} (\bm{\xi}) \mathbf{d}_i = \hat{\mathbf{R}}(\bm{\xi}) \mathbf{d} \label{eq:mid-surf-disp-disretization} \text{ ,}
\end{align}
where 
\begin{align}
    \hat{\mathbf{R}}(\bm{\xi}) = \begin{bmatrix}
        \mathbf{I}^{{sd}} \hat{R}_{1p}(\bm{\xi}) &  \mathbf{I}^{{sd}} \hat{R}_{2p}(\bm{\xi}) & \ldots & \mathbf{I}^{{sd}} \hat{R}_{np}(\bm{\xi})
    \end{bmatrix} \label{eq:shell-nurbs-basis}
\end{align}
is the matrix of NURBS basis function with degree $p$, and $n$ is the number of control points, $\mathbf{I}^{{sd}}$ is the identity matrix in $\mathbb{R}^{{sd}}$ with ${sd}$ as the spatial dimension. We neglect NURBS degree $p$ in the matrix notation for conciseness. The parametric coordinate $\bm{\xi} \in \mathbb{R}^{pd}$, where $pd$ is the parametric dimension. For the isogeometric Kirchhoff--Love shell, $sd=3$ and $pd=2$. $\mathbf{P}_i$ and $\mathbf{d}_i$ are vectors of mid-surface geometry control points and displacements associated with node $i$.  Accordingly, the position vector on the mid-surface shell patch in the deformed configuration given by \eqref{eq:kl-shell-deformed-config} is
\begin{align}
    \mathbf{x}(\bm{\xi}) = \mathbf{X}(\bm{\xi}) + \mathbf{u}(\bm{\xi}) = \hat{\mathbf{R}}(\bm{\xi}) (\mathbf{P}+\mathbf{d}) \text{ .} \label{eq:mid-surf-deformed-disretization}
\end{align}
Substituting \eqref{eq:mid-surf-disp-disretization} and \eqref{eq:mid-surf-deformed-disretization} into \eqref{eq:kl-shell-covariant-basis} and following the procedures \eqref{eq:kl-shell-metric-coeff} -- \eqref{eq:kl-shell-virtual-work}, we can assemble the residual force vector by taking the first derivative of the internal work \eqref{eq:kl-shell-virtual-work} and the stiffness matrix for the second derivative\footnote{We use $\mathrm{d}_{\mathbf{v}}(\cdot)$ and $\partial_{\mathbf{v}}(\cdot)$ to denote the total derivative and partial derivative, respectively, of a function with respect to the discrete variables $\mathbf{v}$. This notation distinguishes from the functional derivative in the continuous setting, denoted as $(\cdot),_{v}$, to avoid confusion.}, respectively,

\begin{align}
    \mathbf{R}_{\text{s}} = \partial_{\mathbf{d}} W_{\text{s}} \qquad \text{and} \qquad \mathbf{K}_{\text{s}} =\partial_{\mathbf{d}} \mathbf{R}_{\text{s}} \label{eq:shell-f-k} \text{ .}
\end{align}
For shell structures with single patch NURBS surface, the displacement increments can be solved by $\mathbf{K}_{\text{s}} \, \Delta\mathbf{d} = -\mathbf{R}_{\text{s}}$. 

\begin{figure}[!htb]\centering
    \includegraphics[width=0.99\textwidth]{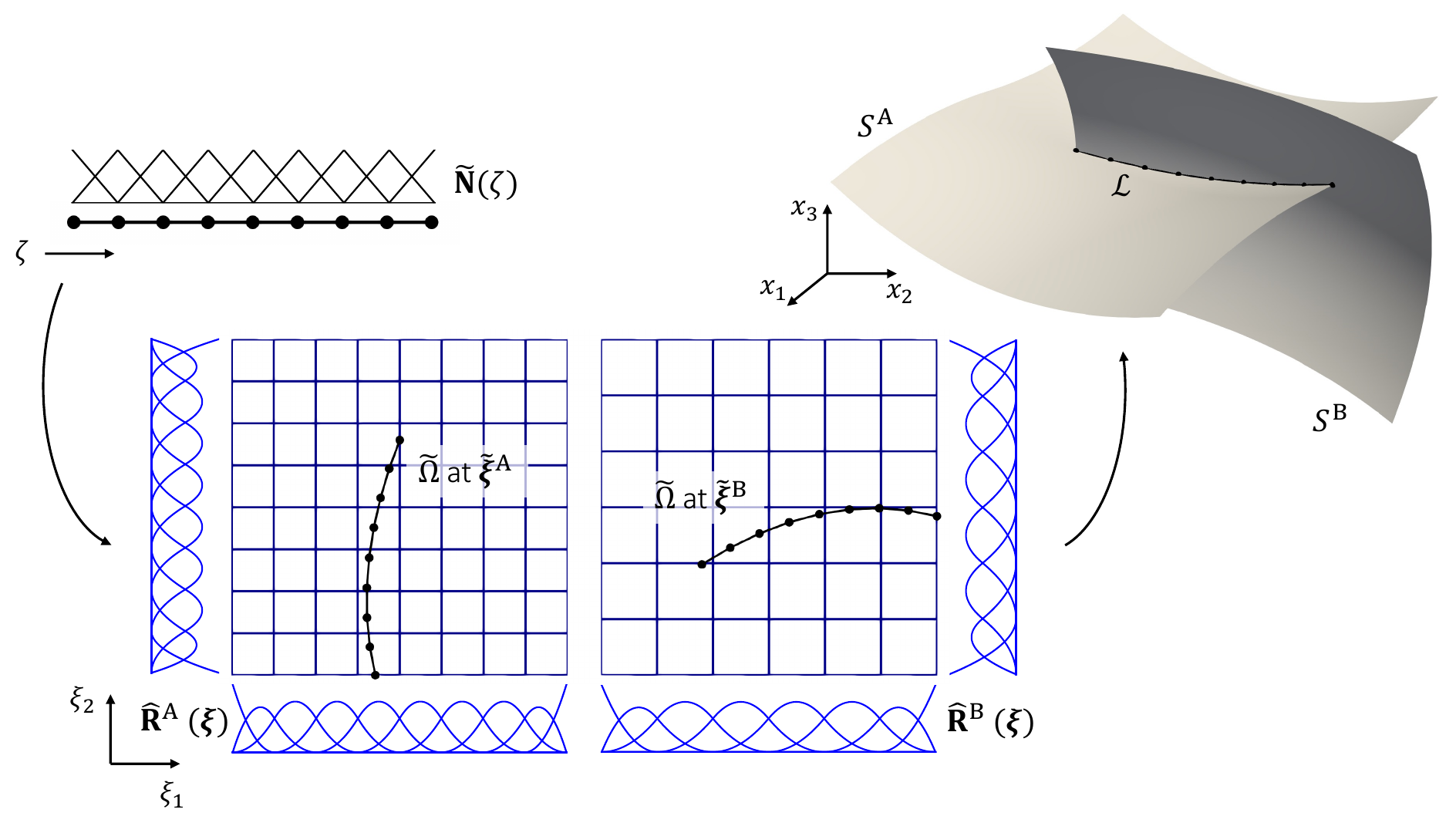}
    \caption{An illustrative example of two shell patches with one intersection. Shell patches $S^\text{A}$ and $S^\text{B}$ are discretized isogeometrically using NURBS basis functions $\hat{\mathbf{R}}^{\text{A}}(\bm{\xi})$ and $\hat{\mathbf{R}}^{\text{B}}(\bm{\xi})$. A topologically 1D quadrature mesh $\tilde{\Omega}$, discretized using linear basis functions $\tilde{\mathbf{N}}(\zeta)$, is created in the parametric space to integrate the penalty energy for shell coupling.}
    \label{fig:two-patch-coupling-example}
\end{figure}

For multi-patch shell structures, contributions of the coupling term outlined in \eqref{eq:penatly-virtual-work} to both membrane and bending stiffness need to be taken into consideration. Using a shell structure with two patches as an example, depicted in Figure \ref{fig:two-patch-coupling-example}, a topologically 1D, geometrically 2D quadrature mesh $\tilde{\Omega}$ \footnote{In this paper, all symbols indicated with $\, \tilde{ }\,$ denote quantities defined on the quadrature mesh of patch intersections.} is constructed in the parameter space to represent the integration domain of the patch intersection. We first move the quadrature mesh to the parametric location of the intersection relative to shell patch $S^{\text{A}}$. The reference geometry and displacements of the patch intersection are obtained by interpolating corresponding functions from $S^{\text{A}}$ to $\tilde{\Omega}$,
\begin{align}
    \tilde{\mathbf{X}}^{\text{A}}(\zeta) = \tilde{\mathbf{N}}(\zeta) \hat{\mathbf{R}}^{\text{A}}(\tilde{\bm{\xi}}^{\text{A}}) \mathbf{P}^{\text{A}} = \tilde{\mathbf{N}}(\zeta) \tilde{\mathbf{P}}^{\text{A}} 
    \qquad \text{and} \qquad 
    \tilde{\mathbf{u}}^{\text{A}}(\zeta) = \tilde{\mathbf{N}}(\zeta) \hat{\mathbf{R}}^{\text{A}}(\tilde{\bm{\xi}}^{\text{A}}) \mathbf{d}^{\text{A}} = \tilde{\mathbf{N}}(\zeta) \tilde{\mathbf{d}}^{\text{A}} \label{eq:intersection-geom-disp} \text{ ,}
\end{align}
where 
\begin{align}
    \tilde{\mathbf{N}}(\zeta) = \begin{bmatrix}
        \mathbf{I}^{sd} \tilde{N}_1(\zeta) & \mathbf{I}^{sd} \tilde{N}_2(\zeta) & \ldots & \mathbf{I}^{sd} \tilde{N}_m(\zeta)
    \end{bmatrix} \label{eq:inersection-physical-basis}
\end{align}
denotes the basis function of the quadrature mesh to approximate quantities in the physical space. Standard liner basis functions are employed for $\tilde{\mathbf{N}}(\zeta)$ in this paper, and $m$ is the number of nodes of the quadrature mesh. $\tilde{\bm{\xi}}^{\text{A}} \in \mathbb{R}^{m\cdot pd}$ refers to the vector of nodal coordinates of the quadrature mesh relative to shell patch $S^{\text{A}}$ with $\tilde{\bm{\xi}}^{\text{A}}_i \in \mathbb{R}^{pd}$. The calculation of $\tilde{\bm{\xi}}^{\text{A}}$ is discussed in Section \ref{subsubsec:implicit-relation-cp-xi}. Additionally, $\zeta$ is the isoparametric coordinate of the quadrature mesh, with $\zeta \in \mathbb{R}^1$ due to $\tilde{\Omega}$ being a topologically 1D mesh. $\hat{\mathbf{R}}^{\text{A}}(\tilde{\bm{\xi}}^{\text{A}}) \in \mathbb{R}^{(m\cdot sd)\times(n\cdot sd)}$ is the interpolation matrix, each row is the evaluation of the NURBS basis function of shell $S^{\text{A}}$ at $\tilde{\bm{\xi}}^{\text{A}}_i$. $\tilde{\mathbf{P}}^{\text{A}}$ and $\tilde{\mathbf{d}}^{\text{A}}$ are vectors of interpolated control points and displacements on the intersection. Substituting \eqref{eq:intersection-geom-disp} into \eqref{eq:kl-shell-deformed-config} and \eqref{eq:kl-shell-covariant-basis}, covariant basis vectors of the mid-surface on the intersection $\mathcal{L}$ are obtained as
\begin{equation}
    \begin{aligned}
    \tilde{\mathbf{A}}^{\text{A}}_{\alpha} &= \tilde{\mathbf{X}}^{\text{A}},_{\xi_\alpha}= \tilde{\mathbf{N}}(\zeta) \hat{\mathbf{R}}^{\text{A}},_{\xi_\alpha}(\tilde{\bm{\xi}}^{\text{A}}) \mathbf{P}^{\text{A}} = \tilde{\mathbf{N}}(\zeta) \tilde{\mathbf{P}}^{\text{A}}_{\xi_{\alpha}} \qquad \text{and} \\[6pt]
    \tilde{\mathbf{a}}^{\text{A}}_{\alpha} &= \tilde{\mathbf{x}}^{\text{A}},_{\xi_\alpha}= \tilde{\mathbf{N}}(\zeta) \hat{\mathbf{R}}^{\text{A}},_{\xi_\alpha}(\tilde{\bm{\xi}}^{\text{A}}) (\mathbf{P}^{\text{A}}+\mathbf{d}^{\text{A}}) = \tilde{\mathbf{N}}(\zeta) (\tilde{\mathbf{P}}^{\text{A}}_{\xi_{\alpha}} + \tilde{\mathbf{d}}^{\text{A}}_{\xi_{\alpha}}) \label{eq:intersection-covariant-vectors} \text{ ,}
    \end{aligned}
\end{equation}
where $\hat{\mathbf{R}}^{\text{A}},_{\xi_\alpha}(\tilde{\bm{\xi}}^{\text{A}})$ is the first order derivative of the interpolation matrix along parametric direction $\xi_{\alpha}$, and $\tilde{\mathbf{P}}^{\text{A}}_{\xi_{\alpha}}$ and $\tilde{\mathbf{d}}^{\text{A}}_{\xi_{\alpha}}$ are interpolated first order derivative of the control points and displacement functions with respect to the parametric coordinates $\tilde{\bm{\xi}}^{\text{A}}$ of intersection $\mathcal{L}$. Plugging \eqref{eq:intersection-covariant-vectors} into \eqref{eq:kl-shell-normal-vectors}, normal vectors of the intersection on shell $S^{\text{A}}$ in the reference and deformed configurations can be computed as $\tilde{\mathbf{A}}^{\text{A}}_3$ and $\tilde{\mathbf{a}}^{\text{A}}_3$. It is notable that \eqref{eq:intersection-covariant-vectors} requires the first order derivatives of the NURBS basis functions, ensuring rotational continuity is preserved at patch intersections.

Tangent vectors of the intersection on both configurations have to be computed before acquiring conormal vectors in \eqref{eq:penatly-virtual-work}, and they are given by
\begin{align}
    \tilde{\mathbf{A}}^{\text{A}}_{t} = \tilde{\mathbf{X}}^{\text{A}},_{\zeta} = \tilde{\mathbf{N}},_{\zeta}(\zeta) \tilde{\mathbf{P}}^{\text{A}}
    \qquad \text{and} \qquad
    \tilde{\mathbf{a}}^{\text{A}}_{t} =  \tilde{\mathbf{x}}^{\text{A}},_{\zeta} = \tilde{\mathbf{N}},_{\zeta}(\zeta) (\tilde{\mathbf{P}}^{\text{A}}+\tilde{\mathbf{u}}^{\text{A}}) \label{eq:intersection-tangent-vectors} \text{ .}
\end{align}
Subsequently, conormal vectors on reference and deformed configurations are defined as
\begin{align}
    \tilde{\mathbf{A}}^{\text{A}}_{n} = \frac{\tilde{\mathbf{A}}^{\text{A}}_{t} \times \tilde{\mathbf{A}}^{\text{A}}_{3}}{\Vert \tilde{\mathbf{A}}^{\text{A}}_{t} \times \tilde{\mathbf{A}}^{\text{A}}_{3} \Vert}
    \qquad \text{and} \qquad
    \tilde{\mathbf{a}}^{\text{A}}_{n} = \frac{\tilde{\mathbf{a}}^{\text{A}}_{t} \times \tilde{\mathbf{a}}^{\text{A}}_{3}}{\Vert \tilde{\mathbf{a}}^{\text{A}}_{t} \times \tilde{\mathbf{a}}^{\text{A}}_{3} \Vert} \label{eq:intersection-conormal-vectors} \text{ .}
\end{align}

Next, we move the quadrature mesh $\tilde{\Omega}$ to the parametric position defined by coordinates $\tilde{\bm{\xi}}^{\text{B}}$ relative to shell patch $S^{\text{B}}$, where the calculation of $\tilde{\bm{\xi}}^{\text{B}}$ is discussed in Section \ref{subsubsec:implicit-relation-cp-xi}. By repeating \eqref{eq:intersection-geom-disp} and \eqref{eq:intersection-covariant-vectors}, we can determine the displacements $\tilde{\mathbf{u}}^{\text{B}}$ and normal vectors $\tilde{\mathbf{A}}^\text{B}_3$ and $\tilde{\mathbf{a}}^\text{B}_3$ of $S^{\text{B}}$ at the intersection $\mathcal{L}$. Substituting these displacements and geometry vectors from the quadrature mesh into \eqref{eq:penatly-virtual-work}, the penalty virtual work $\delta W^{\text{AB}}_{\text{pen}}(\tilde{\mathbf{d}}, \tilde{\mathbf{d}}_{\bm{\xi}}, \tilde{\mathbf{P}}, \tilde{\mathbf{P}}_{\bm{\xi}})$ can be integrated on $\tilde{\Omega}$, where $\tilde{\mathbf{d}}$ and $\tilde{\mathbf{P}}$ are the interpolated displacements and geometric control points for both surfaces. $\tilde{\mathbf{d}}_{\bm{\xi}}$ and $\tilde{\mathbf{P}}_{\bm{\xi}}$ are the associated first order derivatives. Consequently, the residual force vector and stiffness matrix of the coupled shell structure are
\begin{align}
    \mathbf{R}_{} = \begin{bmatrix}
        \mathbf{R}^{\text{A}}_{\text{s}} + \mathbf{R}^{\text{A}}_{\text{pen}}  \\ \mathbf{R}^{\text{B}}_{\text{s}} + \mathbf{R}^{\text{B}}_{\text{pen}}
    \end{bmatrix}
    \qquad \text{and} \qquad
    \mathbf{K}_{} = \begin{bmatrix}
        \mathbf{K}^{\text{A}}_{\text{s}} + \mathbf{K}^{\text{AA}}_{\text{pen}} & \mathbf{K}^{\text{AB}}_{\text{pen}} \\ \mathbf{K}^{\text{BA}}_{\text{pen}} & \mathbf{K}^{\text{B}}_{\text{s}} + \mathbf{K}^{\text{BB}}_{\text{pen}}
    \end{bmatrix}  \label{eq:nonmatching-shell-f-k} \text{ ,}
\end{align}
where components of penalty energy contribution, e.g., $\mathbf{R}^{\text{A}}_{\text{pen}}$ and $\mathbf{K}^{\text{AB}}_{\text{pen}}$, are defined as
\begin{align}
    \mathbf{R}^{\text{A}}_{\text{pen}} &= (\hat{\mathbf{R}}^{\text{A}}(\tilde{\bm{\xi}}^{\text{A}}))^{\mathrm{T}} \, \partial_{\tilde{\mathbf{d}}^{\text{A}}} W^{\text{AB}}_{\text{pen}} + (\hat{\mathbf{R}},_{\bm{\xi}}^{\text{A}}(\tilde{\bm{\xi}}^{\text{A}}))^{\mathrm{T}} \, \partial_{\tilde{\mathbf{d}}_{\bm{\xi}}^{\text{A}}} W^{\text{AB}}_{\text{pen}} \label{eq:nonmatching-shell-f-comp}  \qquad \text{and} \qquad \\[6pt]
    \mathbf{K}^{\text{AB}}_{\text{pen}} &= (\hat{\mathbf{R}}^{\text{B}}(\tilde{\bm{\xi}}^{\text{B}}))^{\mathrm{T}} \, \partial_{\tilde{\mathbf{d}}^{\text{B}}} \mathbf{R}^{\text{A}}_{\text{pen}} + (\hat{\mathbf{R}},_{\bm{\xi}}^{\text{B}}(\tilde{\bm{\xi}}^{\text{B}}))^{\mathrm{T}} \, \partial_{\tilde{\mathbf{d}}_{\bm{\xi}}^{\text{B}}} \mathbf{R}^{\text{A}}_{\text{pen}} \label{eq:nonmatching-shell-k-comp} \text{ .}
\end{align}
And $\hat{\mathbf{R}},_{\bm{\xi}}^{\text{A}}(\tilde{\bm{\xi}}^{\text{A}}) \in \mathbb{R}^{(m \cdot pd \cdot sd)\times (n \cdot sd)}$ is the first order derivative of the interpolation matrix on both parametric directions.

The displacement increments for both spline patches can be solved using the Newton--Raphson method, as expressed by $\mathbf{K}_{} \Delta \mathbf{d} = -\mathbf{R}_{}$. Equation \eqref{eq:nonmatching-shell-k-comp} indicates that $\mathbf{K}^{\text{AB}}_{\text{pen}} = {\mathbf{K}^{\text{BA}}_{\text{pen}}}^{\mathrm{T}}$, enabling the lower triangle blocks in $\mathbf{K}_{}$ to be obtained from the upper triangle counterparts, thereby improving computational efficiency. 
Readers are referred to \cite{zhao2022open} for details about implementation and code framework. A series of benchmark problems in \cite[Section 4]{zhao2022open} have been utilized to verify the accuracy of this method.

\section{Shape optimization of non-matching shells with moving intersections} \label{sec:shape-opt-nonmatching-shells}
Integrating IGA into shell shape optimization presents notable advantages. The direct analysis based on CAD geometries in IGA naturally bridges the gap between the design model and analysis model within the optimization loop without geometric errors. Compared to the classical FEM, IGA-based shape optimization entirely bypasses the process of conforming FE mesh generation, thereby significantly simplifying the workflow due to the absence of FE mesh sensitivity. This section presents the formulations for IGA-based shape optimization, followed by an in-depth discussion of multi-patch shell structures with moving intersections.

\subsection{Shape optimization of isogeometric Kirchhoff--Love shell} \label{subsec:shape-opt-KL-shell}
A general shape optimization problem for an isogeometric shell patch can be formulated as
\begin{equation}
\begin{aligned}
    \underset{\mathbf{P}}{\text{minimize }} &f(\mathbf{P}) \\
    \text{subject to } &\mathbf{g}(\mathbf{P}) \leq \mathbf{0}\\
    &\mathbf{h}(\mathbf{P}) = \mathbf{0} \label{eq:minmize-single-patch} \text{ ,}
\end{aligned}
\end{equation}
where design variable $\mathbf{P}$ are the control points of the shell geometry, $f$ is the objective function, $\mathbf{g}$ and $\mathbf{h}$ are the vector-valued inequality and equality constraints, respectively. We adopt internal energy as the objective function to illustrate the optimization scheme. The internal energy of the Kirchhoff--Love shell is a function of both the control points of geometry $\mathbf{P}$ and displacements $\mathbf{d}$, expressed as $f=W^{\text{int}}_{\text{s}}(\mathbf{P}, \mathbf{d}(\mathbf{P}))$. In this study, a gradient-based optimization algorithm is used due to its benefits in efficiency and suitability to large-scale problems. The total derivative of a single patch shell shape optimization is given by the chain rule
\begin{align}
    \mathrm{d}_{\mathbf{P}}f = \partial_{\mathbf{P}} f + \left(\partial_{\mathbf{d}} f \right)^{\mathrm{T}} \mathrm{d}_{\mathbf{P}}\mathbf{d} \label{eq:expression-dfdp-single-patch-init} \text{ ,}
\end{align}
where the partial derivatives $\partial_{\mathbf{P}} f$ and $\partial_{\mathbf{d}} f $ can be readily calculated with isogeometric discretization in \eqref{eq:mid-surf-disp-disretization}. The total derivative $\mathrm{d}_{\mathbf{P}} \mathbf{d}$ can be determined by the physical constraint of the Kirchhoff--Love shell theory $\mathbf{R}_{\text{s}} (\mathbf{P}, \mathbf{d}) = \mathbf{0}$ for all input $\mathbf{P}$, which implies 
\begin{align}
    \mathrm{d}_{\mathbf{P}}\mathbf{R}_{\text{s}} &= \partial_{\mathbf{P}}\mathbf{R}_{\text{s}} + \partial_{\mathbf{d}}\mathbf{R}_{\text{s}}\, \mathrm{d}_{\mathbf{P}}\mathbf{d} = \mathbf{0} \label{eq:expression-drsdp} \text{ ,} \\
    \mathrm{d}_{\mathbf{P}}\mathbf{d} &= - {(\partial_{\mathbf{d}}\mathbf{R}_{\text{s}})}^{-1}\, \partial_{\mathbf{P}}\mathbf{R}_{\text{s}} = - \mathbf{K}_{\text{s}}^{-1}\,\partial_{\mathbf{P}}\mathbf{R}_{\text{s}} \label{eq:expression-dddp}  \text{ ,}
\end{align}
where $\partial_{\mathbf{P}}\mathbf{R}_{\text{s}}$ represents the partial derivative of the shell residual force vector with respect to geometry control points, and $(\partial_{\mathbf{B}}\mathbf{A})_{ij} = \partial_{\mathbf{B}_j}\mathbf{A}_i$. In the direct method, $\mathrm{d}_{\mathbf{P}}\mathbf{d}$ can be solved with 
\begin{align}
    \mathbf{K}_{\text{s}}\, \mathrm{d}_{\mathbf{P}}\mathbf{d} = -\partial_{\mathbf{P}}\mathbf{R}_{\text{s}} 
    \label{eq:linear-system-dddp} \text{ .}
\end{align} However, the cost of solving \eqref{eq:linear-system-dddp} scales linearly with the number of design variables. The adjoint method is employed to circumvent the increasing expenses of solving the linear systems in \eqref{eq:linear-system-dddp} with a large number of design variables. Substituting \eqref{eq:expression-dddp} into \eqref{eq:expression-dfdp-single-patch}, the total derivative states as
\begin{align}
    \mathrm{d}_{\mathbf{P}}f = \partial_{\mathbf{P}} f - (\partial_{\mathbf{d}} f )^{\mathrm{T}} \mathbf{K}_{\text{s}}^{-1}\,\partial_{\mathbf{P}}\mathbf{R}_{\text{s}} = \partial_{\mathbf{P}} f + (\mathrm{d}_{\mathbf{R}_{\text{s}}}f)^{\mathrm{T}}\,\partial_{\mathbf{P}}\mathbf{R}_{\text{s}}\label{eq:expression-dfdp-single-patch} \text{ ,}
\end{align}
where $\mathrm{d}_{\mathbf{R}_{\text{s}}}f$ can be solved with the following equation \begin{align}
    \mathbf{K}_{\text{s}}^{\mathrm{T}}\, \mathrm{d}_{\mathbf{R}_{\text{s}}}f = - \partial_{\mathbf{d}} f \label{eq:linear-system-dfdrs} \text{ .}
\end{align}
The number of linear solves in \eqref{eq:linear-system-dfdrs} equals the number of model outputs and remains independent of the number of design variables. In practical shape optimization scenarios, the number of design variables typically far exceeds the number of outputs. Therefore, the adjoint method is more advantageous for addressing large-scale optimization problems. By solving the total derivative in \eqref{eq:linear-system-dfdrs} and substituting it into \eqref{eq:expression-dfdp-single-patch}, the shell geometry with minimum internal energy is obtained when the algorithm satisfies the optimality condition.

\subsection{Shape optimization of multi-patch isogeometric Kirchhoff--Love shells} \label{subsec:shape-opt-multi-KL-shell}

Here, we extend the optimization problem \eqref{eq:minmize-single-patch} to encompass multi-patch shell structures, using a two-patch configuration illustrated in Figure \ref{fig:two-patch-coupling-example} to demonstrate the optimization approach. For clarity, we continue to use $\mathbf{P}$ and $\mathbf{d}$ to represent the control points for the geometry and displacements of the non-matching shell. Specifically, we define $\mathbf{P} = \begin{bmatrix}
    {\mathbf{P^{\text{A}}}}^{\mathrm{T}} & {\mathbf{P}^{\text{B}}}^{\mathrm{T}}
\end{bmatrix}^{\mathrm{T}}$ and $\mathbf{d} = \begin{bmatrix}
    {\mathbf{d^{\text{A}}}}^{\mathrm{T}} & {\mathbf{d}^{\text{B}}}^{\mathrm{T}}
\end{bmatrix}^{\mathrm{T}}$. In addition to the change in geometry control points, multi-patch shell structures involve the relative movement between shell patches during shape optimization. To account for this movement, we introduce an additional set of state variables denoted as $\tilde{\bm{\xi}} = \begin{bmatrix}
{\tilde{\bm{\xi}}^{\text{A}}}^{\mathrm{T}} & {\tilde{\bm{\xi}}^{\text{B}}}^{\mathrm{T}}
\end{bmatrix}^{\mathrm{T}}$ as shown in Figure \ref{fig:two-patch-coupling-example}, representing the parametric coordinates of the patch intersections, into the shape optimization process. 

Section \ref{subsec:isogeometric-discretization} indicates that, besides the boundary and load conditions, the displacement field of non-matching shell structures depends not only on the shell geometry but also on the parametric location of patch intersections. This dependence is encapsulated by the shell coupling residual vector $\mathbf{R}_{}(\mathbf{P}, \tilde{\bm{\xi}}, \mathbf{u}) = \mathbf{0}$ introduced in \eqref{eq:nonmatching-shell-f-k}. The total derivative of shape optimization for non-matching shells $\mathrm{d}_{\mathbf{P}} f$ remains the same as given in \eqref{eq:expression-dfdp-single-patch-init}. However, the total derivative $\mathrm{d}_{\mathbf{P}}\mathbf{d}$ is obtained by taking the total derivative of the non-matching residual $\mathbf{R}_{}$,
\begin{align}
    \mathrm{d}_{\mathbf{P}}\mathbf{R}_{} &= \partial_{\mathbf{P}}\mathbf{R}_{} + \partial_{\tilde{\bm{\xi}}} \mathbf{R}_{} \, \mathrm{d}_{\mathbf{P}}{\tilde{\bm{\xi}}} + \partial_{\mathbf{d}} \mathbf{R}_{} \, \mathrm{d}_{\mathbf{P}}{\mathbf{d}} = \mathbf{0} \label{eq:expression-multi-patch-drcoudp} \text{ ,} \\
    \mathrm{d}_{\mathbf{P}}{\mathbf{d}} &= - {(\partial_{\mathbf{d}} \mathbf{R}_{})}^{-1} (\partial_{\mathbf{P}}\mathbf{R}_{} + \partial_{\tilde{\bm{\xi}}} \mathbf{R}_{} \, \mathrm{d}_{\mathbf{P}}{\tilde{\bm{\xi}}}) \label{eq:expression-multi-patch-dddp} \text{ ,}
\end{align}
where $\partial_{\mathbf{d}} \mathbf{R}_{}$ is the stiffness matrix of the non-matching shell, $\partial_{\mathbf{d}} \mathbf{R}_{} = \mathbf{K}_{}$. Similar to the single patch shell, the partial derivative $\partial_{\mathbf{P}}\mathbf{R}_{}$ can be derived from the residual vector of the non-matching shell and has an identical form to $\mathbf{K}_{}$, 
\begin{align}
    \partial_{\mathbf{P}}\mathbf{R}_{} = \begin{bmatrix}
    \partial_{\mathbf{P}^{\text{A}}}\mathbf{R}^\text{A}_{\text{s}} + \partial_{\mathbf{P}^{\text{A}}}\mathbf{R}^\text{A}_{\text{pen}} & \partial_{\mathbf{P}^{\text{B}}}\mathbf{R}^\text{A}_{\text{pen}} \\[6pt] \partial_{\mathbf{P}^{\text{A}}}\mathbf{R}^\text{B}_{\text{pen}} & \partial_{\mathbf{P}^{\text{B}}}\mathbf{R}^\text{B}_{\text{s}} + \partial_{\mathbf{P}^{\text{B}}}\mathbf{R}^\text{B}_{\text{pen}} 
    \end{bmatrix} \label{eq:expression-multi-patch-prcoupp} \text{ ,}
\end{align}
where the blocks related to penalty terms, e.g., $\partial_{\mathbf{P}^{\text{B}}}\mathbf{R}^\text{A}_{\text{pen}}$, can be derived from \eqref{eq:nonmatching-shell-f-comp},
\begin{align}
    \partial_{\mathbf{P}^{\text{B}}}\mathbf{R}^\text{A}_{\text{pen}} = (\hat{\mathbf{R}}^{\text{B}}(\tilde{\bm{\xi}}^{\text{B}}))^{\mathrm{T}} \partial_{\tilde{\mathbf{P}}^{\text{B}}}\mathbf{R}^\text{A}_{\text{pen}} + (\hat{\mathbf{R}},^{\text{B}}_{\bm{\xi}}(\tilde{\bm{\xi}}^{\text{B}}))^{\mathrm{T}} \partial_{\tilde{\mathbf{P}}^{\text{B}}_{\bm{\xi}}}\mathbf{R}^\text{A}_{\text{pen}} \label{eq:expression-multi-patch-drpen-dp} \text{ .}
\end{align}
In contrast to the single patch shell, the non-matching shells require additional derivatives, as indicated in \eqref{eq:expression-multi-patch-dddp}, for shape optimization. The partial derivative $\partial_{\tilde{\bm{\xi}}} \mathbf{R}_{}$ in \eqref{eq:expression-multi-patch-dddp} is crucial for differentiating the movement of the intersection during the shape update of shells. This derivative can be obtained from \eqref{eq:nonmatching-shell-f-comp} since only the penalty terms involve the parametric coordinates of the patch intersection. The derivative is expressed as
\begin{align}
    \partial_{\tilde{\bm{\xi}}} \mathbf{R}_{} = \begin{bmatrix}
    \partial_{\tilde{\bm{\xi}}^{\text{A}}}\mathbf{R}^\text{A}_{\text{pen}} & \partial_{\tilde{\bm{\xi}}^{\text{B}}}\mathbf{R}^\text{A}_{\text{pen}} \\[6pt] \partial_{\tilde{\bm{\xi}}^{\text{A}}}\mathbf{R}^\text{B}_{\text{pen}} & \partial_{\tilde{\bm{\xi}}^{\text{B}}}\mathbf{R}^\text{B}_{\text{pen}} 
    \end{bmatrix} \label{eq:expression-multi-patch-drcoudxi} \text{ ,}
\end{align}
where the detailed derivations of sub-blocks is illustrated in \ref{app:derivative-drpen-dxi} using the chain rule. Upon examination of \eqref{eq:nonmatching-shell-f-comp}, it is apparent that the residual vector of the penalty energy $\mathbf{R}^{\text{A}}_{\text{pen}}$ involves the evaluation of the NURBS basis functions and their first derivatives at parametric coordinates of the patch intersection. Note that \eqref{eq:expression-multi-patch-drcoudxi} necessitates the second-order derivatives for both shell patches,  a condition naturally satisfied by the NURBS functions. Hence, the higher-order continuity in NURBS basis functions not only facilitates direct discretization of the Kirchhoff--Love shell model but also provides a straightforward solution for the relative shell movement in shape optimization problems. This ensures that the optimization process can accurately compute the sensitivities of intersection movements in multi-patch shell structures.

\subsubsection{Implicit relation between shell control points and intersections' parametric coordinates} \label{subsubsec:implicit-relation-cp-xi}

Another derivative that needs to be computed in \eqref{eq:expression-multi-patch-dddp} is the total derivative of parametric coordinates of intersections with respect to shell control points, $\mathrm{d}_{\mathbf{P}}{\tilde{\bm{\xi}}}$. This derivative accounts for the sensitivity of the intersection location $\tilde{\bm{\xi}}$ with respect to the shape changes in shell patches. To obtain the analytical derivatives, we establish a relation between $\tilde{\bm{\xi}}$ and $\mathbf{P}$ through a system of implicit equations. These equations are formulated into a residual vector $\mathbf{R}_{\mathcal{L}}(\mathbf{P}, \tilde{\bm{\xi}})$, which reads 
\begin{align}
    \mathbf{R}_{\mathcal{L}}(\mathbf{P}, \tilde{\bm{\xi}}) = \begin{bmatrix}
        \hat{\mathbf{R}}^{\text{A}}(\tilde{\bm{\xi}}^{\text{A}}_i) \mathbf{P}^{\text{A}} - \hat{\mathbf{R}}^{\text{B}}(\tilde{\bm{\xi}}^{\text{B}}_i) \mathbf{P}^{\text{B}} = \mathbf{0} \\[6pt]
        {L^{\text{A}}_j}^2 - {L^{\text{A}}_{j-1}}^2 = 0 \\[6pt]
        \tilde{\xi}^{\text{A\textbackslash B}}_{k} - 1\backslash 0 = 0 \\[6pt]
        \tilde{\xi}^{\text{A\textbackslash B}}_{l} - 1\backslash 0 = 0
    \end{bmatrix} \quad \begin{array}{l}
         \text{for } i \in \{1,2,\ldots, m\}  \\
         \quad \  j \in \{2,3,\ldots,m-1\}
    \end{array} 
    \label{eq:expression-multi-patch-p-xi-residual} \text{ ,}
\end{align}
where $L^{\text{A}}_j$ is the element length of the quadrature mesh $\tilde{\Omega}$ in physical space defined using the Euclidean distance between two adjacent geometric control points of the quadrature mesh
\begin{align}
    L^{\text{A}}_j = \Vert \hat{\mathbf{R}}^{\text{A}}(\tilde{\bm{\xi}}^{\text{A}}_{j+1}) \mathbf{P}^{\text{A}} - \hat{\mathbf{R}}^{\text{A}}(\tilde{\bm{\xi}}^{\text{A}}_{j}) \mathbf{P}^{\text{A}} \Vert_2 \label{eq:quadrature-mesh-phy-el-size} \text{ .}
\end{align}
The first line of \eqref{eq:expression-multi-patch-p-xi-residual} signifies that the parametric coordinates $\tilde{\bm{\xi}}^{\text{A}}$ and $\tilde{\bm{\xi}}^{\text{B}}$ for node $i$ of $\tilde{\Omega}$ coincide in physical space. This condition ensures the recovery of the same physical intersection curve from the parametric space on sides A and B. The second line of \eqref{eq:expression-multi-patch-p-xi-residual} imposes constraints on the quadrature mesh, requiring equally spaced geometric control points and a uniform physical mesh size. This equation rules out the presence of very small elements in the quadrature mesh. The first two lines of \eqref{eq:expression-multi-patch-p-xi-residual} consist of $4m-2$ equations, while there are $4m$ unknowns in $\tilde{\bm{\xi}}$. 

For an arbitrary intersection between two shell patches subjected to elastic deformation, two discrete points on the interaction parametric coordinates  $\tilde{\bm{\xi}}^{\text{A}}$ or $\tilde{\bm{\xi}}^{\text{B}}$ are located at the edges of the shell surfaces as illustrated in Figure \ref{fig:two-patch-coupling-example}. The last two items in \eqref{eq:expression-multi-patch-p-xi-residual} impose such two additional constraints on interaction kinematics where the two edge coordinates have values of either 1 or 0, depending on their parametric location and are denoted using $1\backslash 0$. The parametric coordinate indices $k$ and $l$ take values of $1,2,2m-1$, or $2m$. These two conditions force the intersection edge points to move along their respective edges during the shape optimization process. Ultimately, the four conditions presented in \eqref{eq:expression-multi-patch-p-xi-residual} guarantee a unique set of intersection parametric coordinates for a given pair of shell surfaces.

With the differentiable residual vector $\mathbf{R}_{\mathcal{L}}$, we can obtain the total derivative $\mathrm{d}_{\mathbf{P}} \tilde{\bm{\xi}}$ using the following expression
\begin{align}
    \mathrm{d}_{\mathbf{P}} \mathbf{R}_{\mathcal{L}} &= \partial_{\mathbf{P}} \mathbf{R}_{\mathcal{L}} + \partial_{\tilde{\bm{\xi}}} \mathbf{R}_{\mathcal{L}} \, \mathrm{d}_\mathbf{P} {\tilde{\bm{\xi}}} = \mathbf{0} \label{eq:expression-multi-patch-drldp} \text{ ,} \\
    \mathrm{d}_\mathbf{P} {\tilde{\bm{\xi}}} &= -(\partial_{\tilde{\bm{\xi}}} \mathbf{R}_{\mathcal{L}})^{-1} \partial_{\mathbf{P}} \mathbf{R}_{\mathcal{L}} \label{eq:expression-multi-patch-dxidp} \text{ ,}
\end{align}
where the partial derivatives $\partial_{\tilde{\bm{\xi}}} \mathbf{R}_{\mathcal{L}}$ and $\partial_{\mathbf{P}} \mathbf{R}_{\mathcal{L}}$ can be readily obtained from \eqref{eq:expression-multi-patch-p-xi-residual}. The derivation for these two partial derivatives is demonstrated in \ref{app:derivative-drl-dxi}.

Substituting \eqref{eq:expression-multi-patch-prcoupp}, \eqref{eq:expression-multi-patch-drcoudxi} and \eqref{eq:expression-multi-patch-dxidp} into \eqref{eq:expression-multi-patch-dddp}, the total derivative of displacements with respect to the geometric control points of the non-matching structures can be obtained. Finally, the total derivative of the non-matching shell shape optimization problem can be computed by substituting \eqref{eq:expression-multi-patch-dddp} into \eqref{eq:expression-dfdp-single-patch-init} to yield
\begin{align}
    \mathrm{d}_{\mathbf{P}} f = \partial_{\mathbf{P}} f - (\partial_{\mathbf{d}} f)^{\mathrm{T}} \mathbf{K}_{}^{-1} \left[\partial_{\mathbf{P}}\mathbf{R}_{} - \partial_{\tilde{\bm{\xi}}} \mathbf{R}_{} \, (\partial_{\tilde{\bm{\xi}}} \mathbf{R}_{\mathcal{L}})^{-1} \partial_{\mathbf{P}} \mathbf{R}_{\mathcal{L}}\right] \label{eq:expression-multi-patch-dfdp} \text{ .}
\end{align}
The $(\partial_{\mathbf{d}} f)^{\mathrm{T}} \mathbf{K}_{}^{-1}$ term can be effectively computed using the adjoint method discussed in Section \ref{subsec:shape-opt-KL-shell}. Depending on the shell discretization and number of points on the intersection quadrature mesh, both the direct method and adjoint method can be considered for calculating $\partial_{\tilde{\bm{\xi}}} \mathbf{R}_{} \, (\partial_{\tilde{\bm{\xi}}} \mathbf{R}_{\mathcal{L}})^{-1} \partial_{\mathbf{P}} \mathbf{R}_{\mathcal{L}}$. 

\begin{figure}[!htb]\centering
    \includegraphics[width=0.99\textwidth]{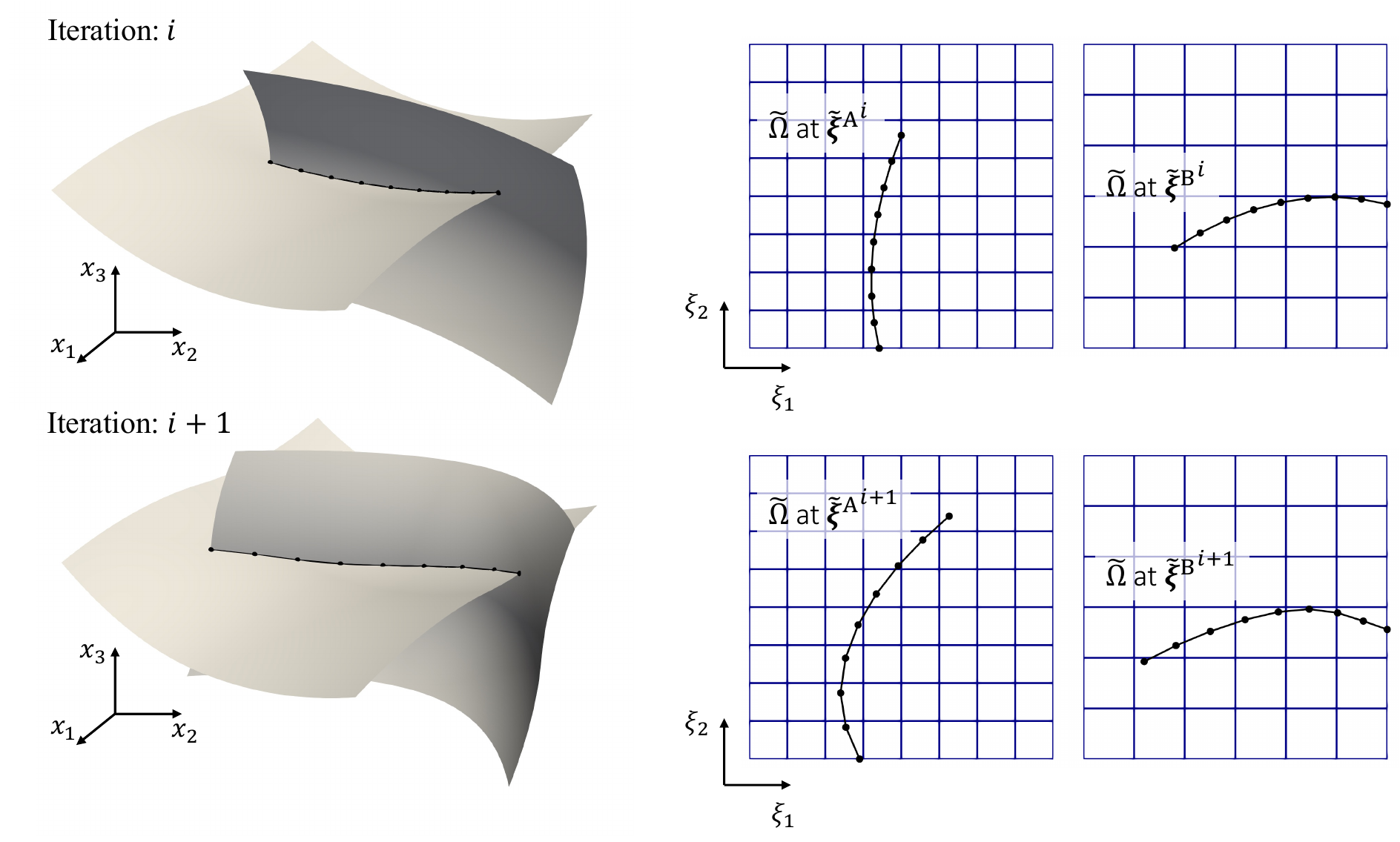}
    \caption{Illustration of shape updates and changes in the relative location for two shell patches during shape optimization. The parametric coordinates of the patch intersection are updated from iteration $i$ to $i+1$ accordingly.}
    \label{fig:shape-opt-scheme}
\end{figure}

By computing the total derivative $\mathrm{d}_{\mathbf{P}}f$ in \eqref{eq:expression-multi-patch-dfdp}, the multi-patch shell structural geometry can be updated using optimization algorithms. A schematic demonstration of the shape update during optimization iterations is depicted in Figure \ref{fig:shape-opt-scheme}.

\section{Implementation}\label{sec:implementation}

The following sections illustrate the implementation details of non-matching shell shape optimization. We adopt the multilevel design concept and present the treatment of various types of intersections. Furthermore, we introduce the dependencies of the open-source Python library used in this paper.

\subsection{Multilevel design for IGA-based optimization}\label{subsec:multilevel-design}

In this paper, we apply the multilevel design concept \cite{nagy2010isogeometric, nagy2013isogeometric, kiendl2014isogeometric} to create a flexible design space. The optimizer modifies only the shape of shell structures with coarse discretizations, referred to as the design model, by adjusting the coordinates of their control points. Meanwhile, a refined geometry, named the analysis model, is used for accurate analysis of the structural response after shape modifications. Specifically, for CAD geometries defined using NURBS basis functions, order elevation ($p$-refinement), knot refinement ($h$-refinement), and the combination of these two methods ($k$-refinement) can be employed to produce finer models while preserving the original geometry. This capability in IGA is particularly beneficial for shape optimization problems as it allows the dimension of the design space to be chosen independently from the dimension of the analysis model. Notably, this approach does not introduce geometric errors into the optimization problem.

\begin{figure}[!htb]\centering
    \includegraphics[width=0.99\textwidth]{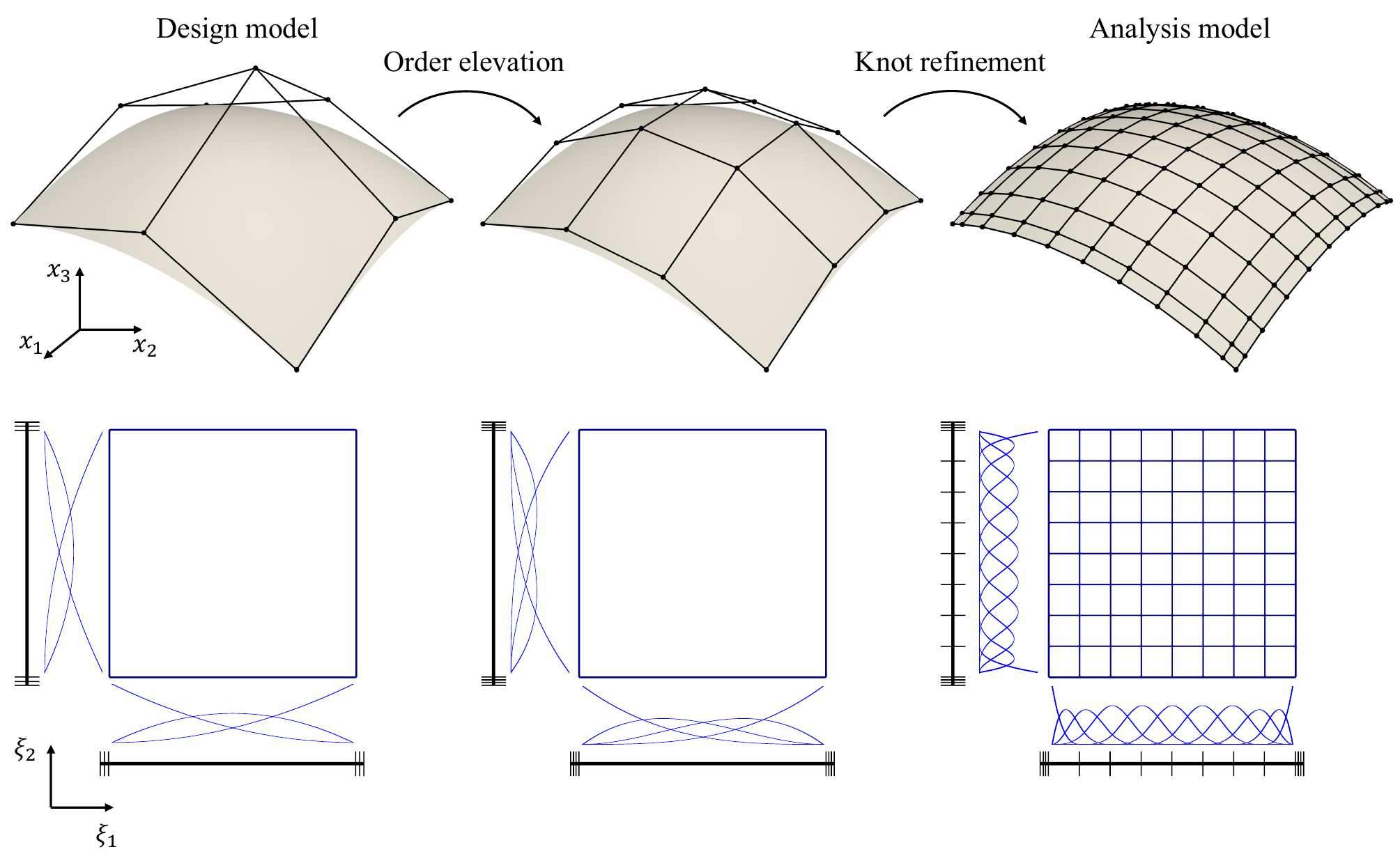}
    \caption{Multilevel design approach for shape optimization problems. The coarse design model is employed to update the shape of the geometry, while the refined analysis model is used for structural analysis. Both the design model and the analysis model represent the same geometry.}
    \label{fig:multilevel-design-concept}
\end{figure}

Figure \ref{fig:multilevel-design-concept} presents an example of the multilevel design approach for a single patch shell. A quadratic surface with coarse discretization is defined by control points $\mathbf{P}^{\text{DM}}$ and NURBS basis functions $\hat{\mathbf{R}}^{\text{DM}}(\bm{\xi})$, which is characterized by a knot vector $\Xi^{\text{DM}} = [[0,0,0,1,1,1], [0,0,0,1,1,1]]$. We first increase the order of NURBS basis functions from quadratic to cubic by adding two extra knots on each side in $\Xi^{\text{DM}}$. The cubic NURBS basis functions $\hat{\mathbf{R}}^{\text{OE}}(\bm{\xi})$ are determined by knot vector $\Xi^{\text{OE}} = [[0,0,0,0,1,1,1,1], [0,0,0,0,1,1,1,1]]$. Consequently, the control points of the surface after order elevation, $\mathbf{P}^{\text{OE}}$, are defined as
\begin{align}
    \hat{\mathbf{R}}^{\text{DM}}(\bm{\xi}) \mathbf{P}^{\text{DM}} = \hat{\mathbf{R}}^{\text{OE}}(\bm{\xi}) \mathbf{P}^{\text{OE}} \text{ .}
\end{align}
With $\hat{\mathbf{R}}^{\text{OE}}$ and $\Xi^{\text{OE}}$, we can insert a sequence of new knots $[0.125, 0.25, 0.375, 0.5, 0.625, 0.75, 0.875]$ into $\Xi^{\text{OE}}$ on both parametric directions to obtain a $h$-refined model with basis functions $\hat{\mathbf{R}}^{\text{KR}}(\bm{\xi})$ characterized by knots vector $\Xi^{\text{KR}}$. The geometric control points of the $h$-refined model are calculated by
\begin{align}
    \hat{\mathbf{R}}^{\text{OE}}(\bm{\xi}) \mathbf{P}^{\text{OE}} = \hat{\mathbf{R}}^{\text{KR}}(\bm{\xi}) \mathbf{P}^{\text{KR}} \text{ .}
\end{align}
We can achieve significantly more accurate analysis results by employing $\hat{\mathbf{R}}^{\text{KR}}(\bm{\xi})$ and $\mathbf{P}^{\text{KR}}$ in the analysis compared to the design model without altering its geometry. Meanwhile, the design model, which has much fewer degrees of freedom (DoFs) compared to the analysis model, allows improved convergence for optimization problems. Design engineers also have the flexibility to define the dimension of the design space by selecting the initial knot vector. It is noted that the continuity in $\hat{\mathbf{R}}^{\text{KR}}$ is increased by one from $\hat{\mathbf{R}}^{\text{DM}}$ through the combination of order elevation and knot refinement, which is advantageous for problems with higher-order governing equations. Commonly used algorithms for the implementation of these refinement strategies are introduced in the NURBS book \cite[Chapter 5]{piegl2012nurbs}.

The multilevel design approach can be readily extended to shape optimization with non-matching shell structures, where the differentiation for the movement of patch intersections during the optimization process is discussed in Section \ref{subsec:shape-opt-multi-KL-shell}.

\subsection{Intersection types in shape optimization} \label{subsec:intersection-types}

\begin{figure}[!htb]
    \centering
    \begin{subfigure}[t]{0.32\textwidth}
        \centering
        \includegraphics[width=\textwidth]{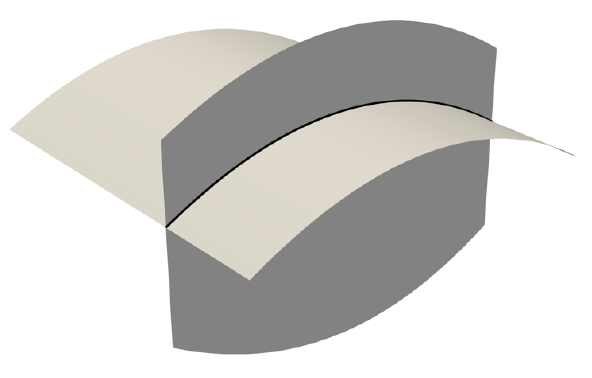}
        \caption{Intersection type: interior--interior}
        \label{subfig:int-type-surf-surf}
    \end{subfigure}
    \hfill
    \begin{subfigure}[t]{0.32\textwidth}
        \centering
        \includegraphics[width=\textwidth]{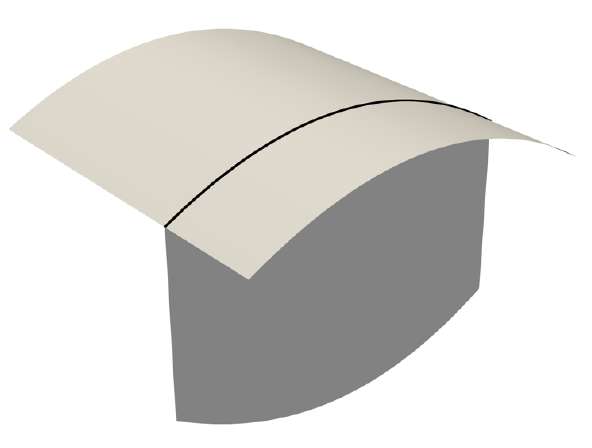}
        \caption{Intersection type: interior--edge}
        \label{subfig:int-type-surf-edge}
    \end{subfigure}
    \hfill
    \begin{subfigure}[t]{0.32\textwidth}
        \centering
        \includegraphics[width=\textwidth]{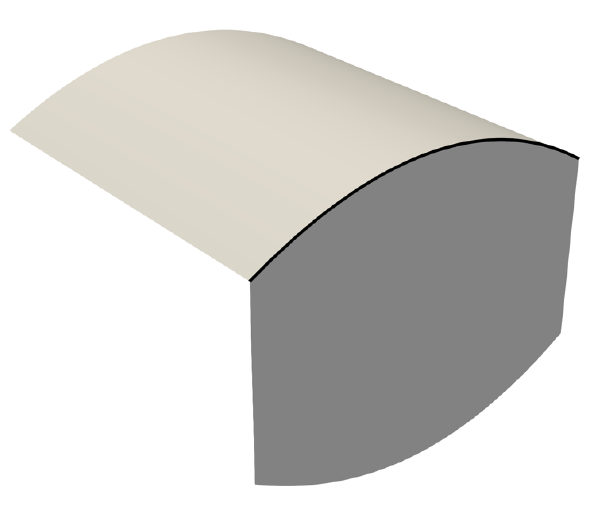}
        \caption{Intersection type: edge--edge}
        \label{subfig:int-type-edge-edge}
    \end{subfigure}
    \caption{Types of shell patch intersection in shape optimization problems.}
    \label{fig:surf-int-types}
\end{figure}

Without considering extreme cases such as singular points and singular curves, there are typically three types of intersections between two tensor-product NURBS patches, as shown in Figure \ref{fig:surf-int-types}. The first type, named interior--interior intersection, is depicted in Figure \ref{subfig:int-type-surf-surf}. For the interior--interior intersections, we assume that two shell patches can move independently of each other without any other constraints imposed. The second type, termed as interior--edge intersection and illustrated in Figure \ref{subfig:int-type-surf-edge}, occurs when the edge of one shell patch intersects the interior of the other shell patch, forming a T-junction structure. During the optimization process, the intersection is allowed to move while maintaining the T-junction. Therefore, an additional constraint is necessary to fulfill this requirement. For the third intersection type, as shown in Figure \ref{subfig:int-type-edge-edge}, the edges from two separate patches join together and no relative movement between the two patches is allowed. In this intersection topology, the optimization framework enforces the conditions that the relative location of the intersection remains fixed and the two shell patches are always connected at their edges. While these intersection topologies do not represent all possible geometries, they are effective within our targeted applications, particularly in the context of aircraft wing design.

\begin{figure}[!htb]\centering
    \includegraphics[width=0.6\textwidth]{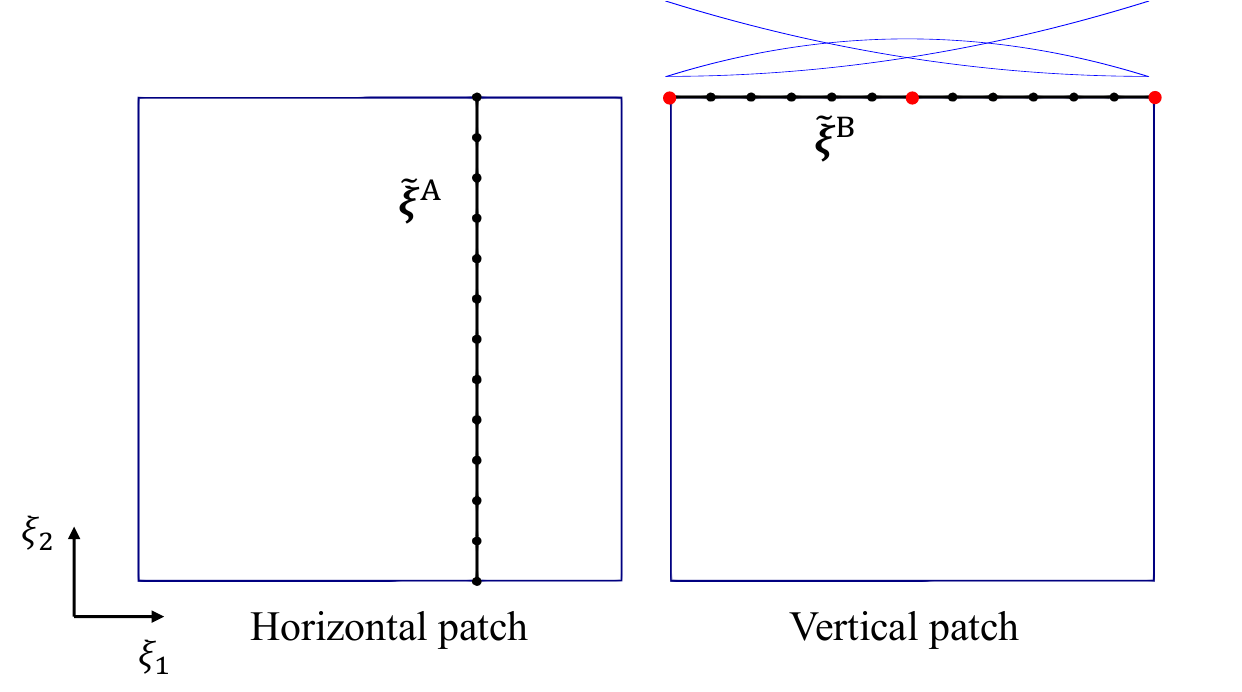}
    \caption{Parametric configuration of two shell patches with an interior--edge intersection.}
    \label{fig:edge-constraint}
\end{figure}

For the interior--edge type of intersections, a linear constraint is applied to the parametric coordinates of the intersection to retain the T-junction. Figure \ref{fig:edge-constraint} depicts the associated parametric configuration and the intersection's quadrature mesh of Figure \ref{subfig:int-type-surf-edge}. To preserve the T-junction, the quadrature mesh related to the vertical patch needs to stay on the top edge. Assuming the parametric domain of the vertical patch is a unit square and the lower-left corner is at $(0,0)$, the constraint is applied as $\tilde{\xi}^\text{B}_{i2}=1$ for $i \in \{1,2,\ldots,m\}$. In the example shown in Figure \ref{fig:edge-constraint}, the intersecting edge of the vertical patch is only defined by three DoFs, leading to an over-constrained system since $m>3$. Therefore, we select three points, highlighted in red, in the quadrature mesh to enforce the T-junction constraint. The support of each NURBS basis function at the intersecting edge needs to contain at least one selected point to uniquely define the edge. It is noted that the edge alignment of the vertical and horizontal patches is imposed only at the selected points to avoid an over-constrained condition. Given the potential for high-order polynomial intersections between two shell patches, the determined curve is considered as an approximated intersecting edge within the design space.

In cases where shell patches form edge-edge intersections, the coordinates of the quadrature mesh are assumed to remain unchanged throughout the optimization process. If the optimization problem incorporates the shape of these shell patches, we employ the FFD-based method as proposed in \cite[Section 4]{zhao2024automated} to ensure the connectivity between shell patches. The shell patches with edge--edge intersections are embedded within a trivariate B-spline block where the shape of shells is updated through the change of the 3D B-spline block. Meanwhile, parametric coordinates of intersections between shell patches in different B-spline blocks are allowed to move. This strategy is employed in the tube optimization benchmark problem in Section \ref{subsec:benchmark-tube}. Conversely, if the shell patches with edge--edge intersections are not considered in the optimization problem, their control points can be fixed without any updates.

\subsection{Software elements for open-source implementation} \label{subsec:software-elements}
The shape optimization Python library is developed leveraging a suite of open-source code packages. It employs the Python interface of OpenCASCADE, PythonOCC \cite{paviot2018pythonocc}, to import the CAD geometry in IGES or STEP formats into the optimization process. Meanwhile, the surface--surface intersection approximation functionality in PythonOCC is utilized to determine the parametric coordinates of intersections, which serve as the initial guess for \eqref{eq:expression-multi-patch-p-xi-residual}. For automated structural analysis of CAD geometries consisting of non-matching isogeometric Kirchhoff--Love shells, the FEniCS \cite{logg2012automated}-based library PENGoLINS \cite{zhao2022open} is employed. The Lagrange polynomial basis functions in the finite element code of FEniCS are changed to NURBS basis functions through the extraction technique \cite{BSEH11, SBVSH11, Schillinger2016, fromm2023interpolation}. The Lagrange extraction is implemented in tIGAr \cite{Kamensky2019}, while the low-level assembly subroutines in FEniCS are reused in the analysis framework.

FEniCS makes use of advanced code generation and computer algebra to automate analytical Gateaux derivative computation, allowing for large-scale gradient-based optimization. Partial derivatives in \eqref{eq:expression-multi-patch-dfdp} are encapsulated into individual components, and they are modularized through OpenMDAO \cite{gray2019openmdao} to manage the adjoint method of total derivative calculation. For solving the optimization problem, the SLSQP optimizer \cite{kraft1988software} is used for simple benchmark examples. The SNOPT optimizer \cite{gill2005snopt}, renowned for its efficiency in nonlinear problems where gradient evaluations are computationally intensively, is employed for complicated problems. The sparse sequential quadratic programming (SQP) algorithm is used in the SNOPT optimizer. The source code of the shape optimization framework is publicly available on the GitHub repository GOLDFISH \cite{goldfish-code}, where demonstrations presented in Sections \ref{sec:benchmark} and \ref{sec:applications} can be reproduced.

\subsection{Optimization scheme} \label{subsec:optimization-scheme}
\begin{figure}[!htb]\centering
    \includegraphics[width=0.35\textwidth]{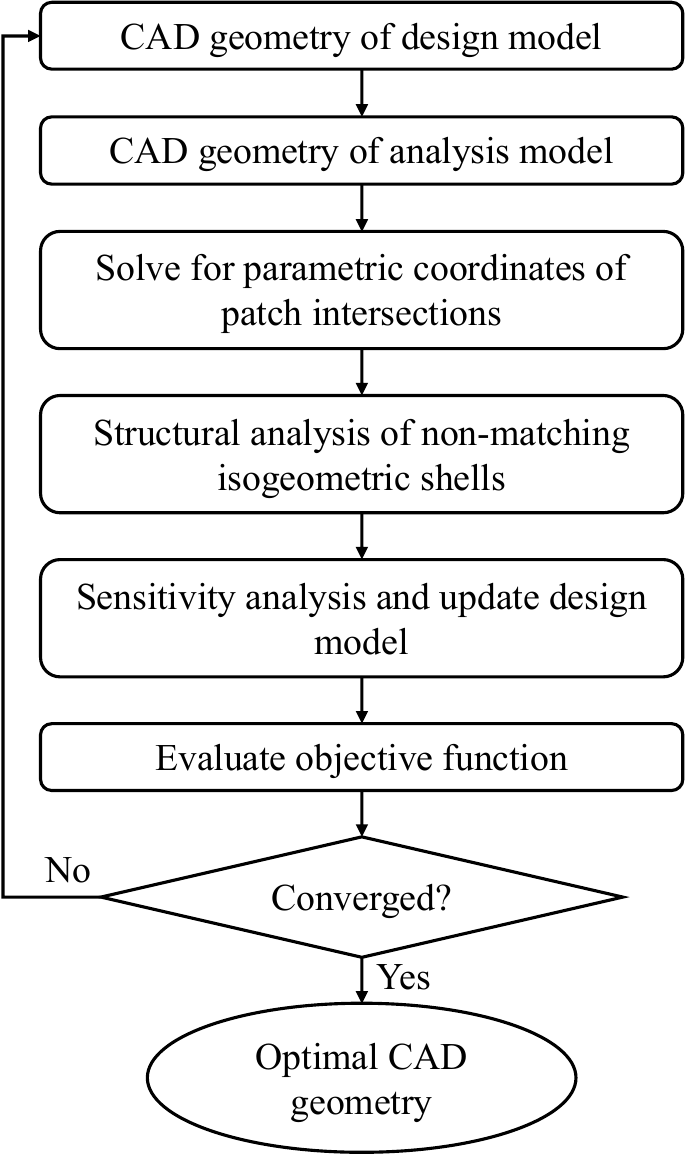}
    \caption{Workflow of the IGA-based shape optimization for non-matching shell structures with moving intersections.}
    \label{fig:optimization-workflow}
\end{figure}

With the aforementioned implementation details and code dependencies, the workflow of shape optimization for non-matching shells is outlined in Figure \ref{fig:optimization-workflow}. The optimization workflow entirely bypasses the FE mesh generation for the CAD geometry. Shape modifications are directly applied to the coarse design model, and the structural response of the updated geometry is evaluated using the refined analysis model. As such, the dimension of the design space can be significantly reduced. As discussed in Section \ref{subsec:multilevel-design}, the geometry preservation properties of NURBS surface refinement methods ensure that no geometric errors are introduced from the design model to the analysis model, which is difficult to achieve in traditional FEM. Consequently, this optimization approach guarantees both accurate geometry representation and analysis results.

\section{Benchmark problems}\label{sec:benchmark}
In this section, we present results based on a set of shape optimization problems to validate the effectiveness of the proposed optimization scheme. The multilevel design approach is employed in the T-beam example, while the FFD-based method, which maintains edge--edge intersections, is tested in the tube problem.

\subsection{T-beam under distributed load}\label{subsec:bechmark-tbeam}
Two types of T-beam geometry are demonstrated to verify the accuracy of the shape optimization approach. The T-beam geometry in Section \ref{subsubsec:straight-tbeam} has a flat top surface, while the top surface in Section \ref{subsubsec:curved-tbeam} is curved to test the proposed approach's ability to preserve the T-junction in curved structure in the optimization process. In both demonstrations, the T-beam is subjected to a downward distributed pressure and is fixed at the rear end. 

\subsubsection{Flat T-beam}\label{subsubsec:straight-tbeam}
\begin{figure}[!htb]
    \centering
    \begin{subfigure}[t]{0.4\textwidth}
        \centering
        \includegraphics[width=\textwidth]{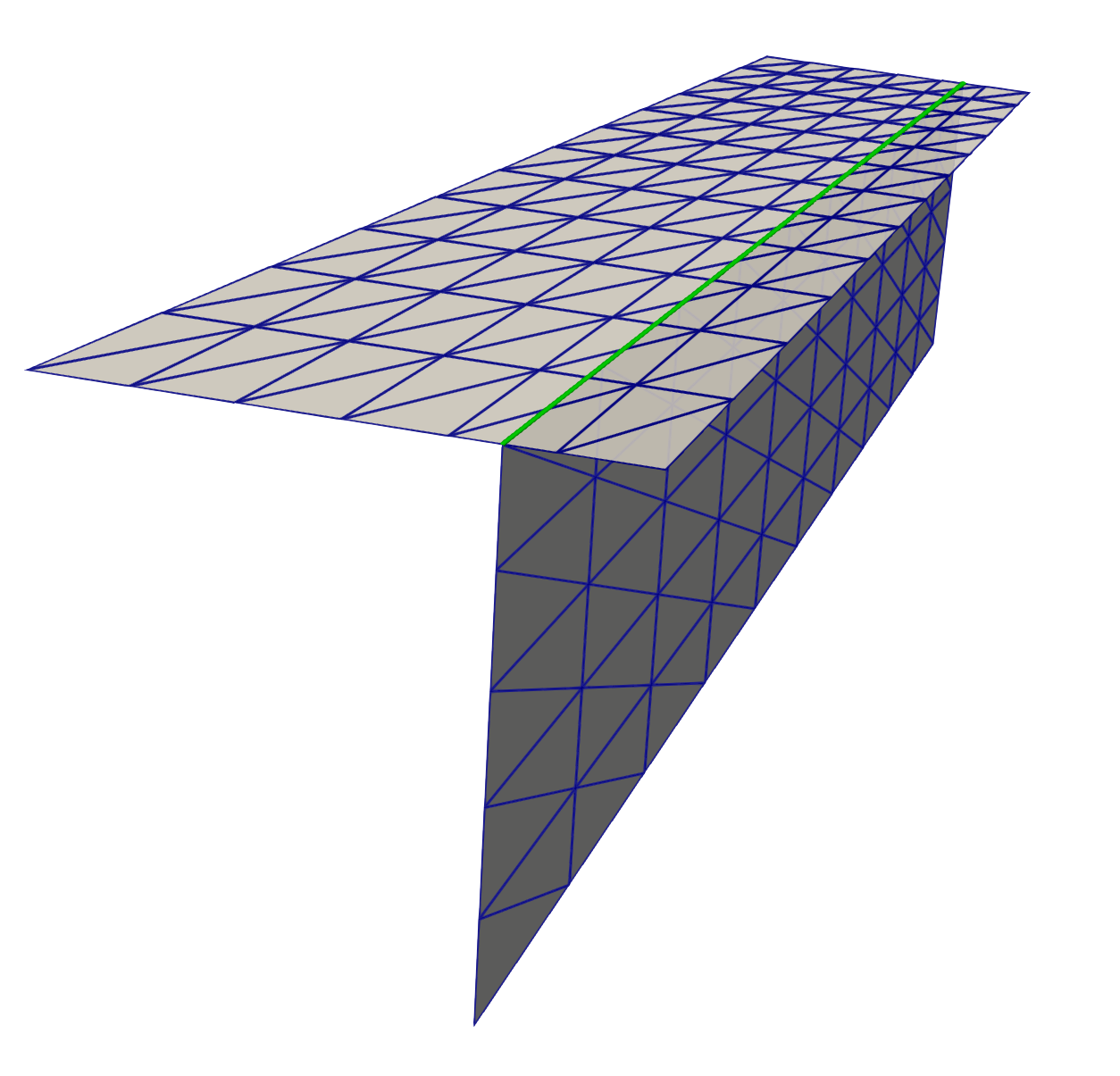}
        \caption{}
        \label{subfig:tbeam-init-geom-flat}
    \end{subfigure}
    \hfill
    \begin{subfigure}[t]{0.55\textwidth}
        \centering
        \includegraphics[width=\textwidth]{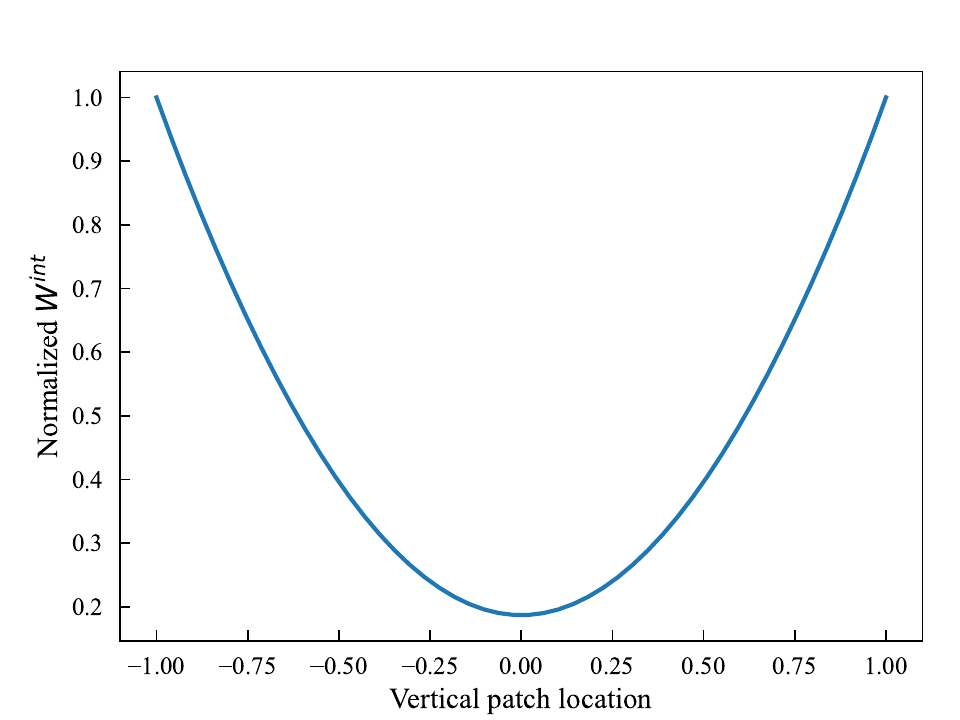}
        \caption{}
        \label{subfig:tbeam-wint-test}
    \end{subfigure}
    \caption{(a) The isogeometrically discretized T-beam geometry with a flat top surface in the initial configuration. The movable intersection is highlighted with a green line. (b) The T-beam's internal energy depends on the location of the vertical surface. The minimal internal energy occurs when the vertical surface is located at the center of the top surface.}
    \label{fig:T-beam-flat}
\end{figure}

For the first benchmark problem, we consider a T-beam geometry composed of two patches, a top surface and a vertical surface. In the optimization process, both surfaces remain flat, with dimensions of 2 m in width and 10 m in length for each patch. The thickness of both shell patches is set as 0.1 m. In the initial design, the top surface ranges from -1 m to 1 m in the horizontal direction, while the top edge of the vertical patch is located at 0.5 m horizontal location of the horizontal patch. The isogeometrically discretized\footnote{Due to technical limitations within FEniCS, the interpolation matrix described in \eqref{eq:intersection-geom-disp} can only be constructed with triangular meshes in the current implementation. While all numerical examples are discretized using triangular meshes, the solutions are still approximated using NURBS basis functions.} analysis model using cubic NURBS basis functions is shown in Figure \ref{subfig:tbeam-init-geom-flat}, where the interior--edge intersection is indicated with a green line. Material properties, Young's modulus $E=10^7\ \rm{Pa}$ and Poisson's ratio $\nu=0$, are used in the analysis, and the uniformly distributed load has a magnitude of $P=1\ \rm{Pa}$.

\begin{figure}[!htb]\centering
    \includegraphics[width=0.9\textwidth]{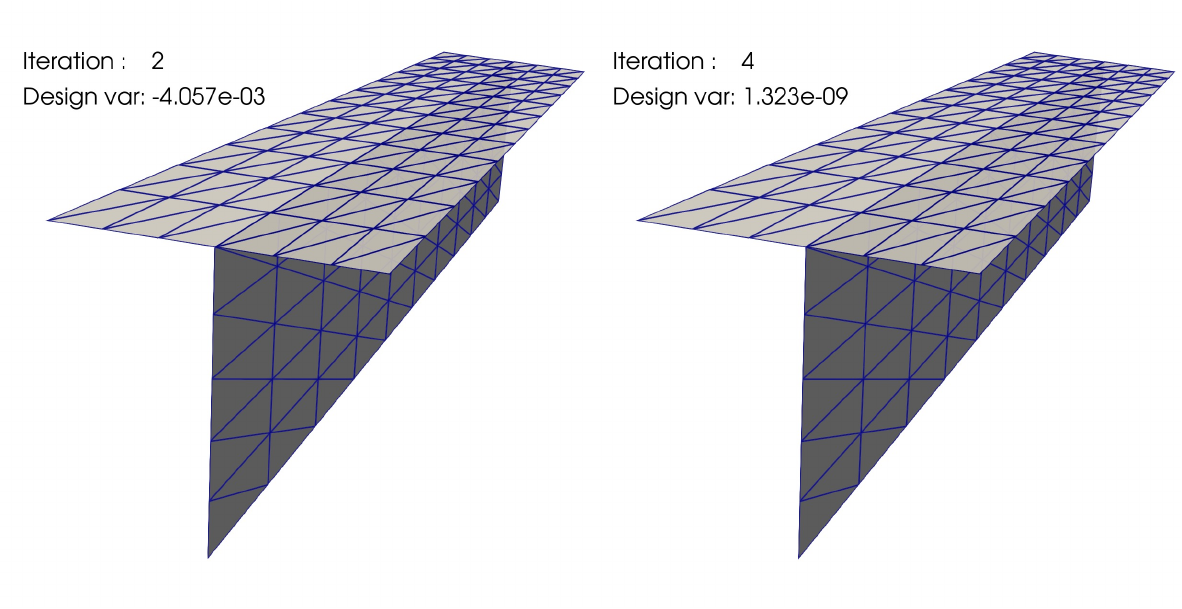}
    \caption{Optimization history of the T-beam with a flat top surface. The optimizer with a tolerance of $10^{-15}$ terminates after 4 iterations.}
    \label{fig:tbeam-flat-shopt-history}
\end{figure}

In this benchmark problem, we aim to minimize the internal energy of the T-beam by adjusting the position of the vertical patch. Thus, only one design variable is considered in this problem. The relation between the internal energy of the T-beam and the location of the vertical patch is illustrated in Figure \ref{subfig:tbeam-wint-test}. The lowest normalized internal energy, with a value of 0.18719, corresponds to the vertical patch positioned at the center of the top patch. Since the movement of the vertical patch is restricted to the horizontal direction, the requirement for the maintenance of the T-junction is automatically satisfied, and the volume of the T-beam remains constant. The only required constraints in this problem are the limits for the coordinate of the vertical patch, which ranges from -1 m to 1 m. The SLSQP optimizer is adopted for this problem with a tolerance set as $10^{-15}$. Due to the simplicity of this benchmark example, the optimizer converges to the optimal location rapidly and terminates successfully with 4 iterations. Two snapshots of the shape update history are demonstrated in Figure \ref{fig:tbeam-flat-shopt-history}. In the converged geometry, the vertical patch has a horizontal coordinate of $1.323\times 10^{-9}$, closely matching the theoretical optimal solution of 0 with a negligible difference. The normalized internal energy of the converged solution has a value of 0.18721, which shows good agreement with the expected value.

\subsubsection{Curved T-beam}\label{subsubsec:curved-tbeam}
\begin{figure}[!htb]\centering
    \includegraphics[width=0.4\textwidth]{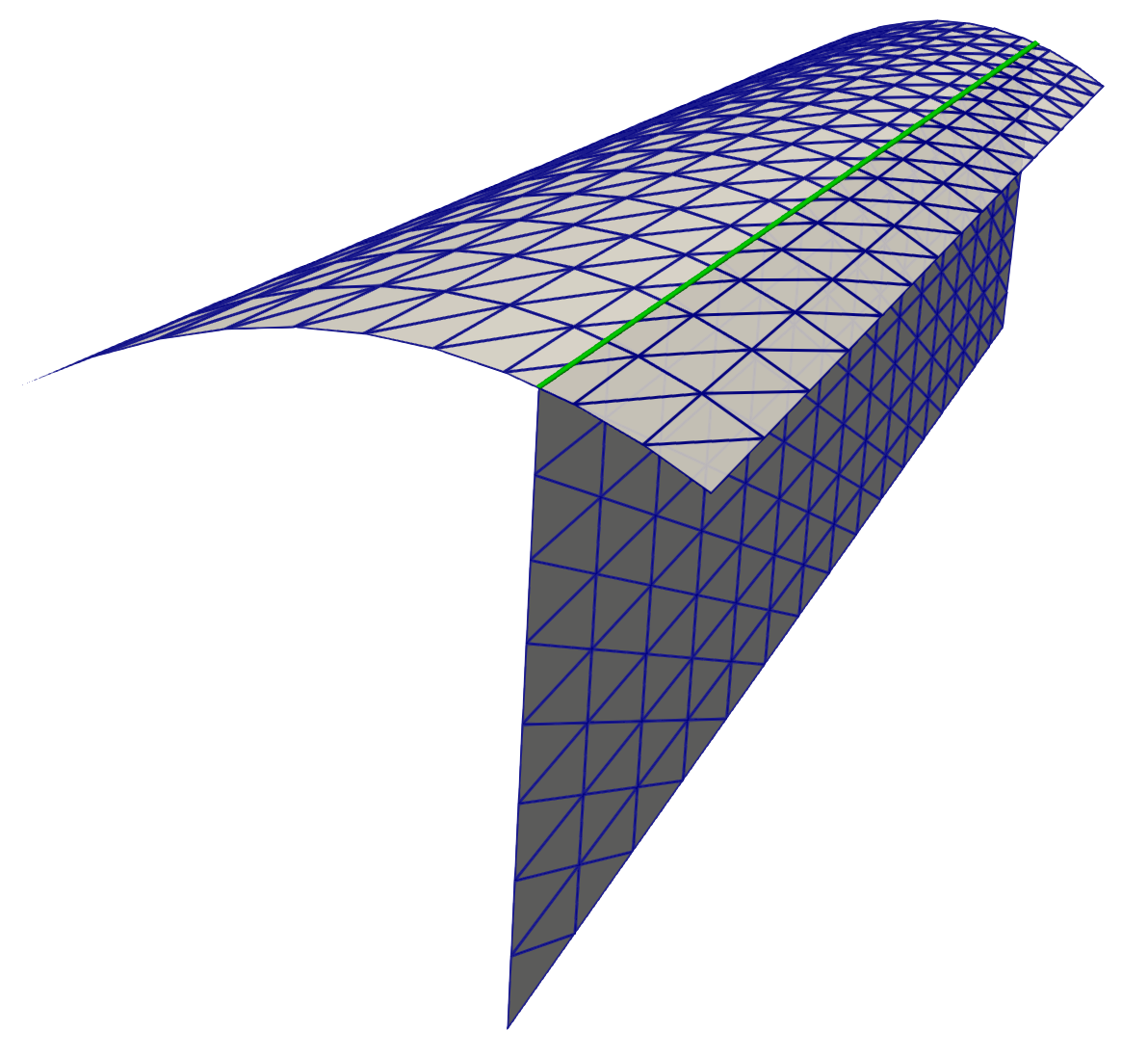}
    \caption{Initial configuration of a T-beam geometry with a curved top surface, where the green line indicates the initial location of the intersection.}
    \label{fig:init-tbeam-curved-geom}
\end{figure}

For this purpose, a T-beam CAD geometry with a curved top surface is generated, and the associated analysis model discretized with cubic NURBS basis functions is shown in Figure \ref{fig:init-tbeam-curved-geom}. The top surface ranges horizontally from -1 m to 1 m and vertically from 0 m to 0.3 m. The vertical surface is located at 0.5 m horizontal location in the initial configuration, where the intersection is marked by a green line. In this benchmark problem, the dimensions of the design space are increased. In the design model, we employ a cubic NURBS curve with a knot vector of $[0, 0, 0, 0, 1, 1, 1, 1]$ to define the horizontal position of the vertical patch, alongside a linear NURBS curve with a knot vector $[0,0,1,1]$ for its vertical location. The vertical patch remains straight in the axial direction during the optimization. Thus, this problem involves four horizontal and two vertical design variables. The same material parameters and objective functions as in Section \ref{subsubsec:straight-tbeam} are used. A constraint ensuring that the top edge of the vertical surface remains attached to the top surface during the optimization is introduced by fixing the parametric coordinate of the quadrature mesh with respect to the vertical patch to 1.0 in the $\xi_2$ direction. Additionally, a volume constraint is imposed on the vertical surface to ensure a constant volume.

\begin{figure}[!htb]\centering
    \includegraphics[width=0.99\textwidth]{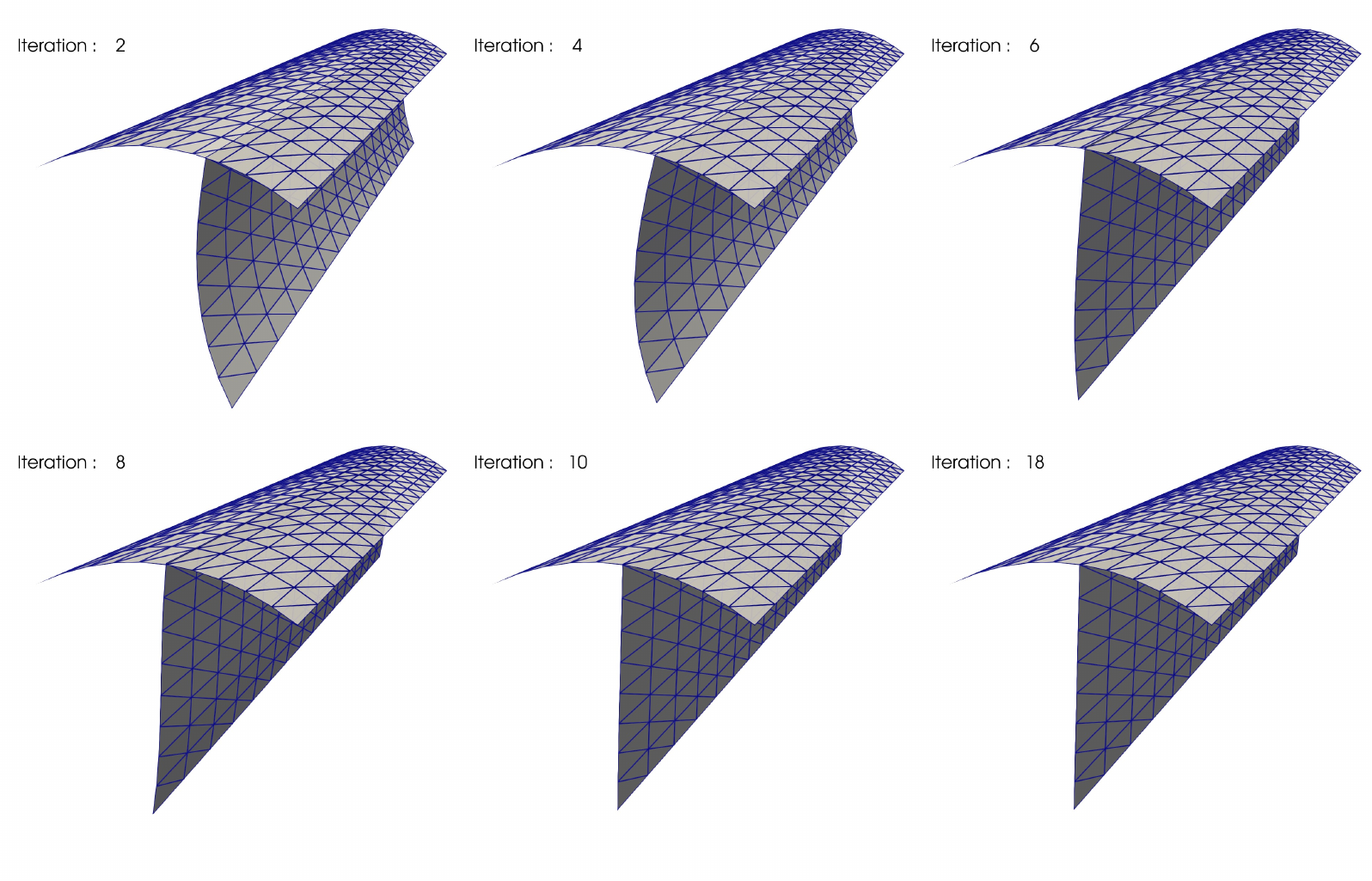}
    \caption{Snapshots of the shape optimization history of the T-beam featuring a curved top surface. The SLSQP optimizer requires 18 iterations to converge to the tolerance of $10^{-12}$.}
    \label{fig:tbeam-curved-shopt-history}
\end{figure}

We continued to employ the SLSQP optimizer with a tolerance of $10^{-12}$. It takes 18 iterations for the optimizer to converge to the specified tolerance. A series of representative optimization snapshots of this benchmark problem is shown in Figure \ref{fig:tbeam-curved-shopt-history}, which demonstrates that the top edge of the vertical surface remains adhered to the top surface due to the implementation of the T-junction preservation constraint. Despite the increased dimension of the design space allowing for potential bending of the vertical patch, it eventually converges to a flat surface in the optimal configuration to minimize internal energy. The coordinates of the four horizontal control points in the optimized design are $[-2.9780\times10^{-9}, -1.0846\times10^{-9}, -2.2580\times10^{-9}, -2.7100\times10^{-9}]$, which correspond to the flat vertical surface at the center of the top surface with sufficiently small errors. Meanwhile, the coordinates of the two vertical control points in the optimal design, $[0.3, -1.7]$, exhibit errors within the machine precision, indicating the volume of the vertical surface remains constant. Accordingly, the vertical coordinate 0.3 denotes that the top edge of the vertical surface precisely lies in the middle of the top surface, maintaining the T-junction connection.

To demonstrate the effectiveness of the proposed shape optimization approach for multi-patch shell structures that incorporate a moving intersection, we test it against two T-beam benchmark problems. Both benchmarks converge to the optimal shapes with sufficiently small errors. During the optimization process, relative movement between the surface patches is achieved using analytical derivatives calculated from the adjoint method, as discussed in Section \ref{subsec:shape-opt-multi-KL-shell}. Additionally, the T-junction is accurately preserved through a linear constraint applied to the parametric coordinates of the intersection's quadrature mesh.

\subsection{Tube with follower pressure}\label{subsec:benchmark-tube}
\begin{figure}[!htb]
    \centering
    \begin{subfigure}[t]{0.4\textwidth}
        \centering
        \includegraphics[width=\textwidth]{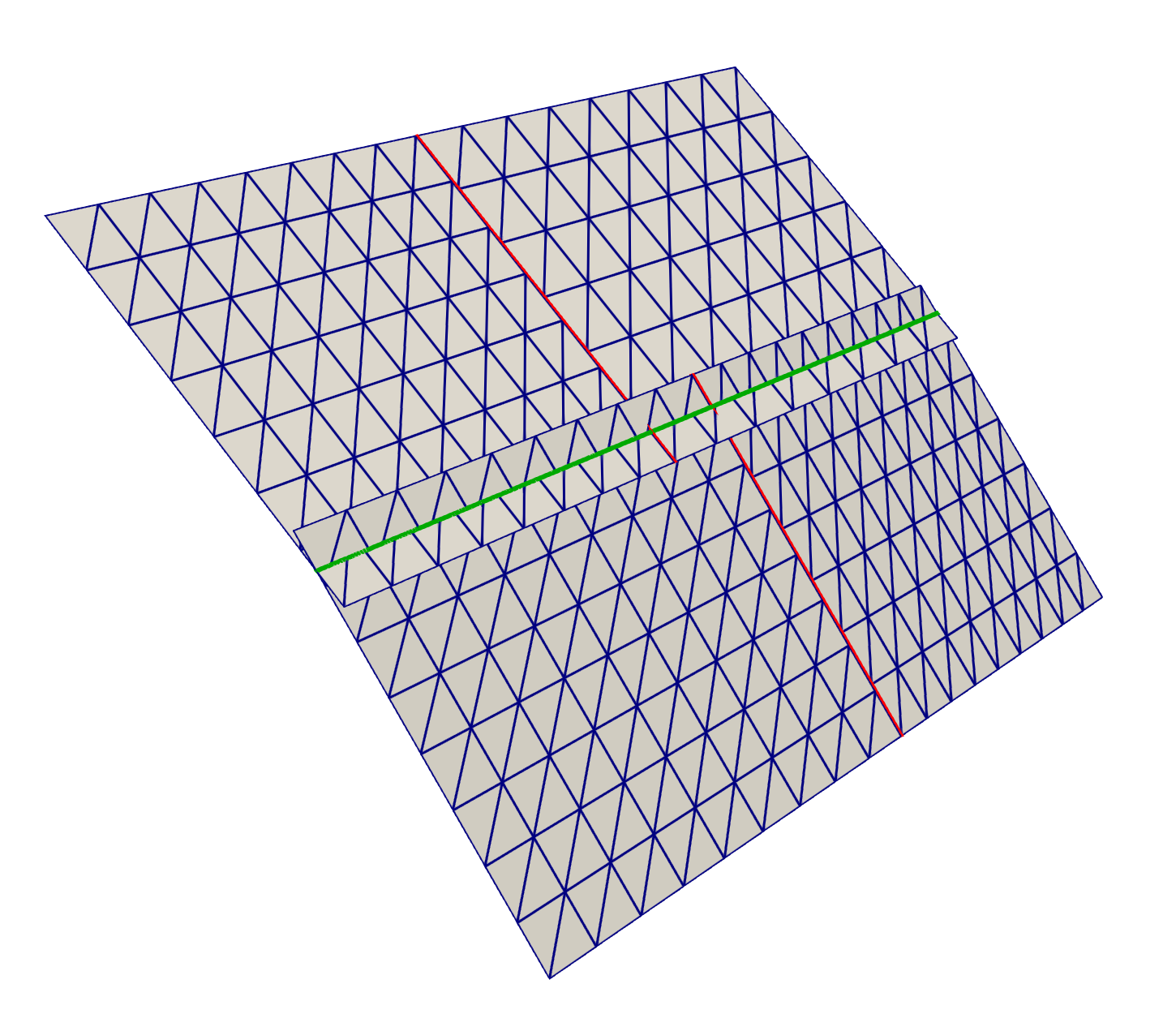}
        \caption{}
        \label{subfig:tube-geom}
    \end{subfigure}
    \hfill
    \begin{subfigure}[t]{0.49\textwidth}
        \centering
        \includegraphics[width=\textwidth]{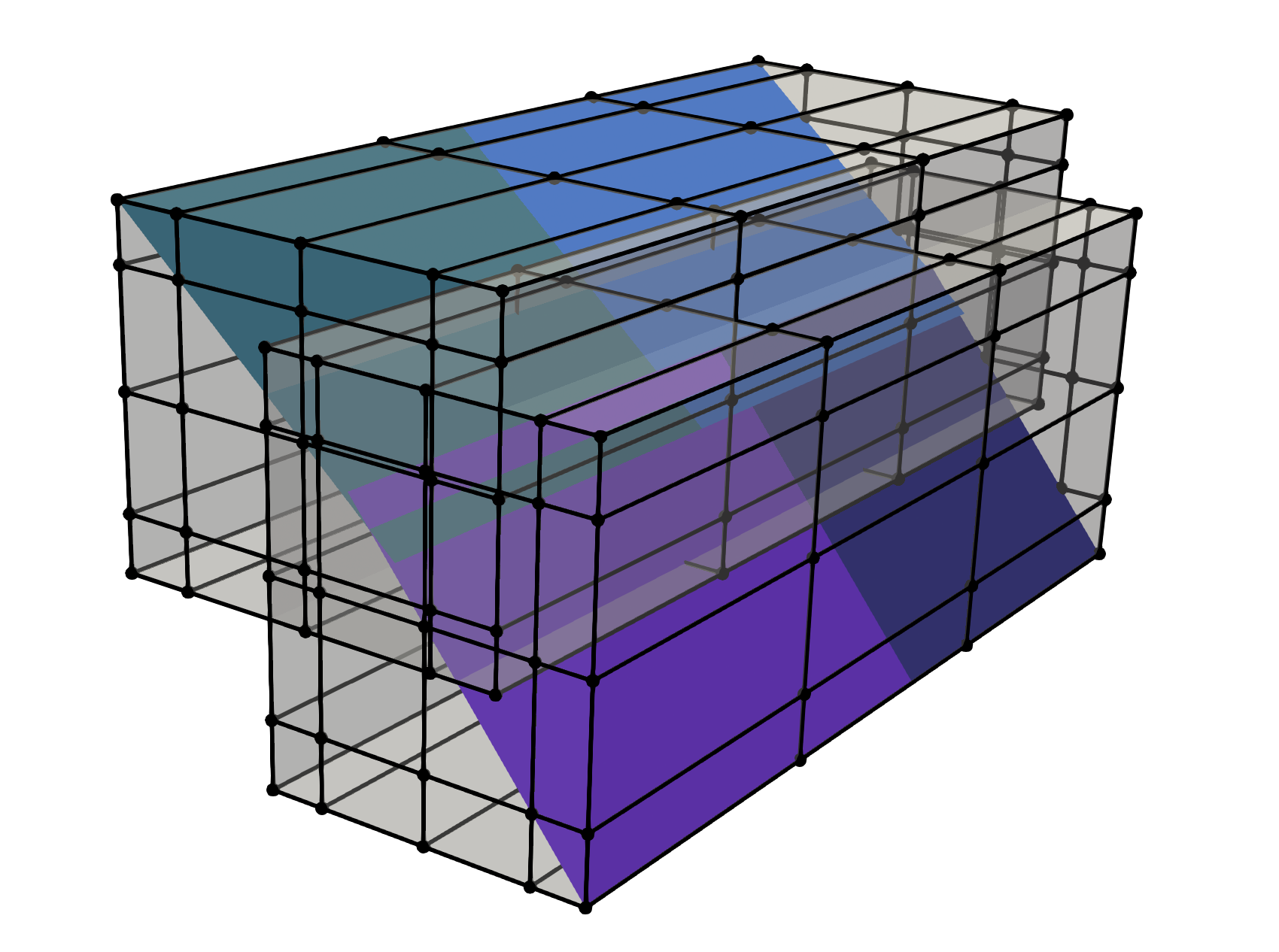}
        \caption{}
        \label{subfig:tube-ffd-init}
    \end{subfigure}
    \caption{(a) A quarter of the initial tube geometry consists of four non-matching cubic NURBS patches. The three interior--interior intersections are indicated with green lines, and two red lines mark the edge--edge intersections. (b) Initial configuration of the tube geometry and FFD blocks. Each set of surface patches with edge--edge intersections is embedded in one 3D B-spline FFD block to preserve the edge--edge intersection, while the interior--interior intersections between different FFD blocks are allowed to move during the shape optimization process.}
    \label{fig:tube-shopt-setup}
\end{figure}

In this section, we investigate the shape optimization of a tube subjected to an outward-facing follower unit pressure on its inner surface. We model a quarter of the tube geometry using four separately parametrized surfaces, with the initial quarter tube geometry depicted in Figure \ref{subfig:tube-geom}. Symmetric boundary conditions are applied to represent the full tube. The initial tube geometry features five intersections in total, two edge--edge intersections, highlighted with red lines, and three interior--interior intersections, marked with green lines. As discussed in Section \ref{subsec:intersection-types}, we assume that the edge--edge intersections remain unchanged due to lack of relative movement, and their intersection type does not alter throughout the optimization process. On the other hand, interior--interior intersections can be moved during the shape optimization, allowing for the search of optimal intersection locations. Consequently, the upper two shell patches can move relative to the lower two patches, and the relative locations within each pair are maintained.

In this benchmark problem, we employ the FFD-based shape modification strategy, incorporating the Lagrange extraction technique \cite{Schillinger2016}, as introduced in \cite{zhao2024automated} for automated preservation of edge--edge intersections in the upper and lower shell patch pairs. The setup of the B-spline blocks in the initial configuration are demonstrated in Figure \ref{subfig:tube-ffd-init}, where the four shell patches are distinguished by different colors. The initial quarter tube geometry ranges from 0 m to 1 m in both vertical and horizontal directions, and from 0 m to 2 m in the axial direction. Each shell patch pair is embedded in a trivariate B-spline block, with shape updates of shell patches achieved by adjusting the control points of the B-spline blocks. Due to the continuous shape modification inside the B-spline block, the edge--edge intersections are maintained. Moreover, relative movement is allowed between the distinct FFD B-spline blocks assigned to the upper and lower pairs. 

\begin{figure}[!htb]
    \centering
    \begin{subfigure}[t]{0.4\textwidth}
        \centering
        \includegraphics[width=\textwidth]{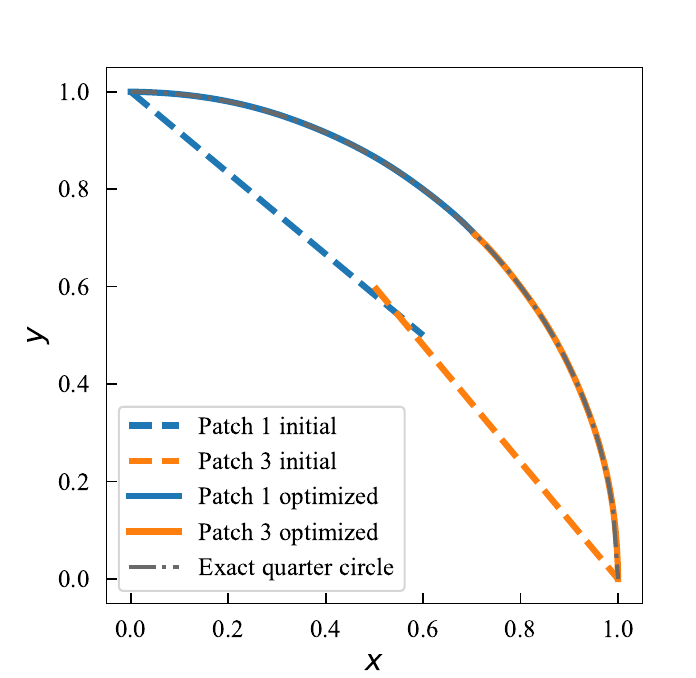}
        \caption{}
        \label{subfig:tube-slice}
    \end{subfigure}
    \hfill
    \begin{subfigure}[t]{0.49\textwidth}
        \centering
        \includegraphics[width=\textwidth]{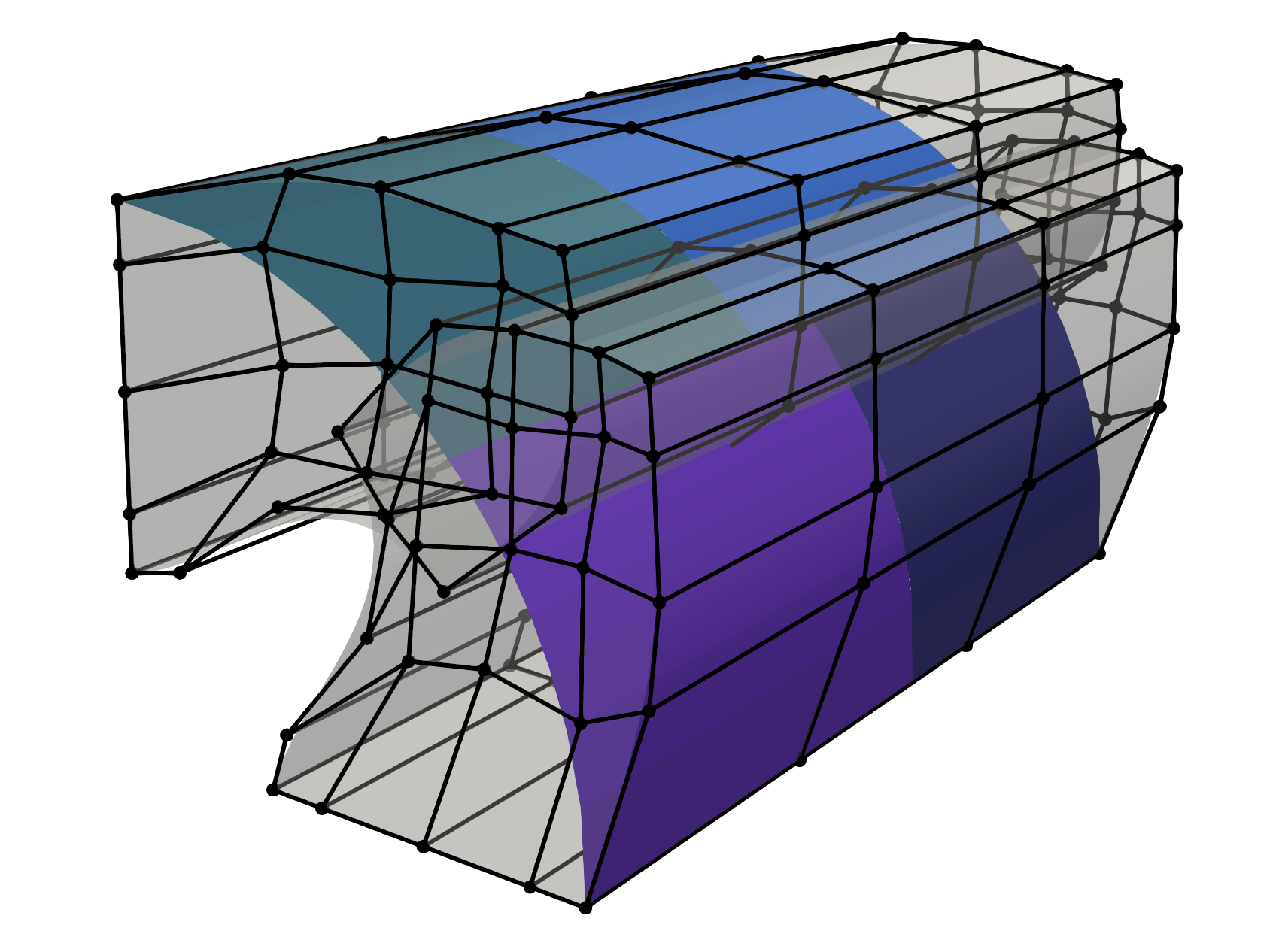}
        \caption{}
        \label{subfig:tube-ffd-final}
    \end{subfigure}
    \caption{(a) Cross-sectional view of the tube geometry in the initial and optimized configurations. (b) Optimized configuration of the tube geometry and FFD blocks.}
    \label{fig:tube-shopt-final}
\end{figure}

In the structural analysis, we use a Young's modulus of $10^{9}\ \rm{Pa}$ and Poisson's ratio of $0$ for the material properties of shell patches, each with a thickness of $0.01\ \rm{m}$. Control points of each FFD block are aligned in the axial direction to ensure the tube remains straight, leading to the assignment of the control points in the first layer of the FFD blocks along the axial direction as design variables. In sum, there are 50 design variables in total, 25 for each FFD block. Meanwhile, the left edge of the upper FFD block and the lower edge of the lower FFD block are fixed to ensure constant positioning of the symmetric edges in the tube geometry. We employ the SNOPT optimizer with a tolerance of $10^{-2}$, requiring 142 iterations to achieve convergence. Figure \ref{fig:tube-shopt-final} displays a sequence of snapshots for the optimization process, with the red curve indicating the cross-section of an exact circular tube. The circular shape represents the theoretical optimal shape that minimizes internal energy under the given follower pressure load conditions. The optimization snapshots demonstrate a gradual transition of the initial tube toward the expected circular tube. Notably, the upper pair of shell patches move freely relative to the lower pair during the optimization iterations. As the optimization progresses, the intersections between these shell pairs shift from interior positions in the initial configuration to the edges in the final configuration, eventually achieving the optimal design.

\begin{figure}[!htb]\centering
    \includegraphics[width=0.99\textwidth]{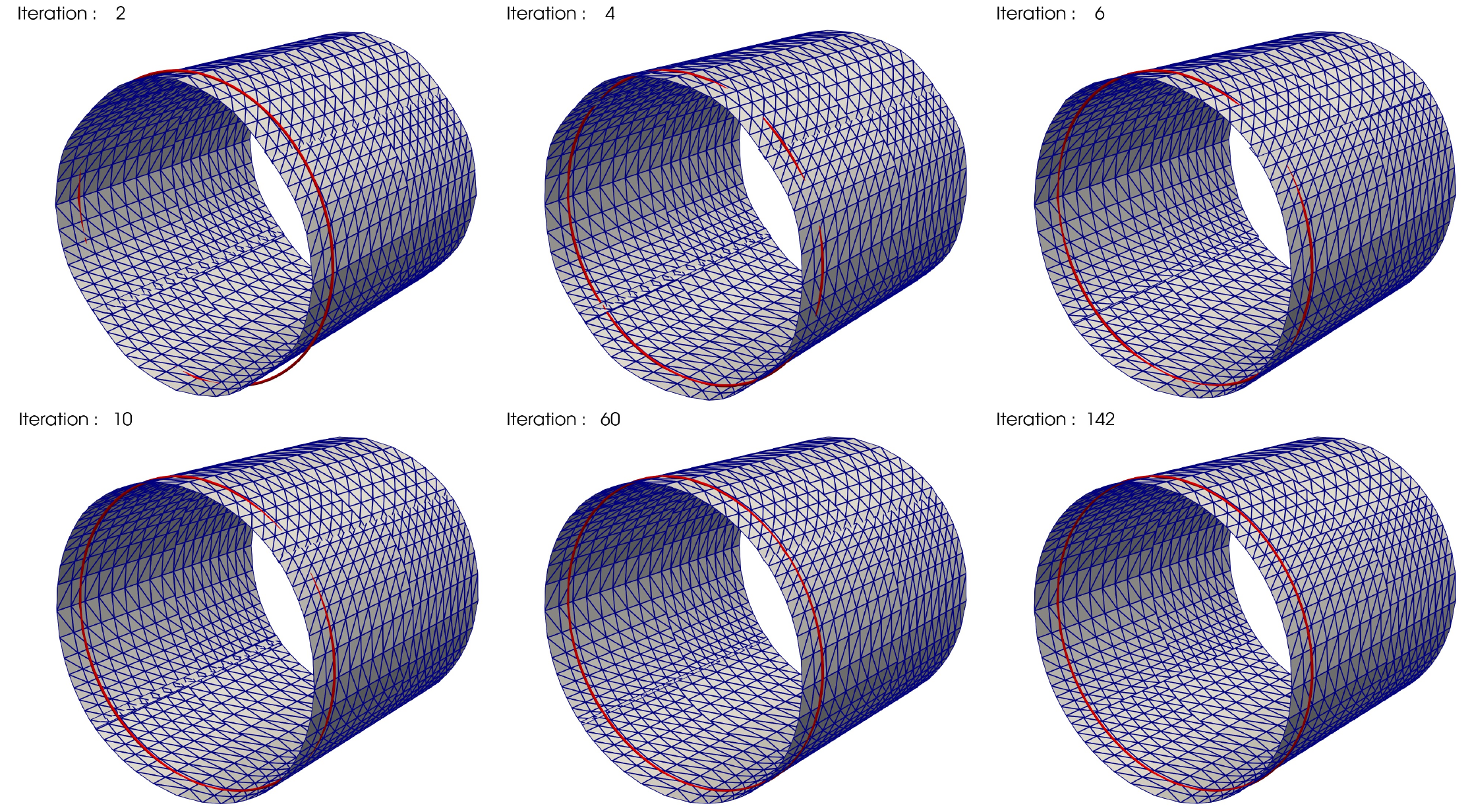}
    \caption{Representative snapshots of the shape optimization history for the tube geometry with interior--interior intersections. The interior--interior intersections converge to edge--edge intersections in the optimized design to minimize the internal energy of the tube.}
    \label{fig:tube-shopt-history}
\end{figure}

A comparison of the cross-sectional view of the tube in initial and optimized configurations is shown in Figure \ref{subfig:tube-slice}. In the optimized configuration, the cross-section of the tube geometry aligns closely with a perfect quarter circular arc, demonstrating the accuracy of the optimization approach. Additionally, the shape of the FFD blocks and associated control points in the final state are displayed in Figure \ref{subfig:tube-ffd-final}. This tube benchmark problem highlights the capability of the optimization approach for handling intersections of the interior--interior type. This approach allows the associated intersecting shell patches to move independently, subject to a constraint guaranteeing the existence of the intersection during the shape optimization process.

\section{Applications to aircraft wings}\label{sec:applications}
\begin{figure}[!htb]\centering
    \includegraphics[width=0.75\textwidth]{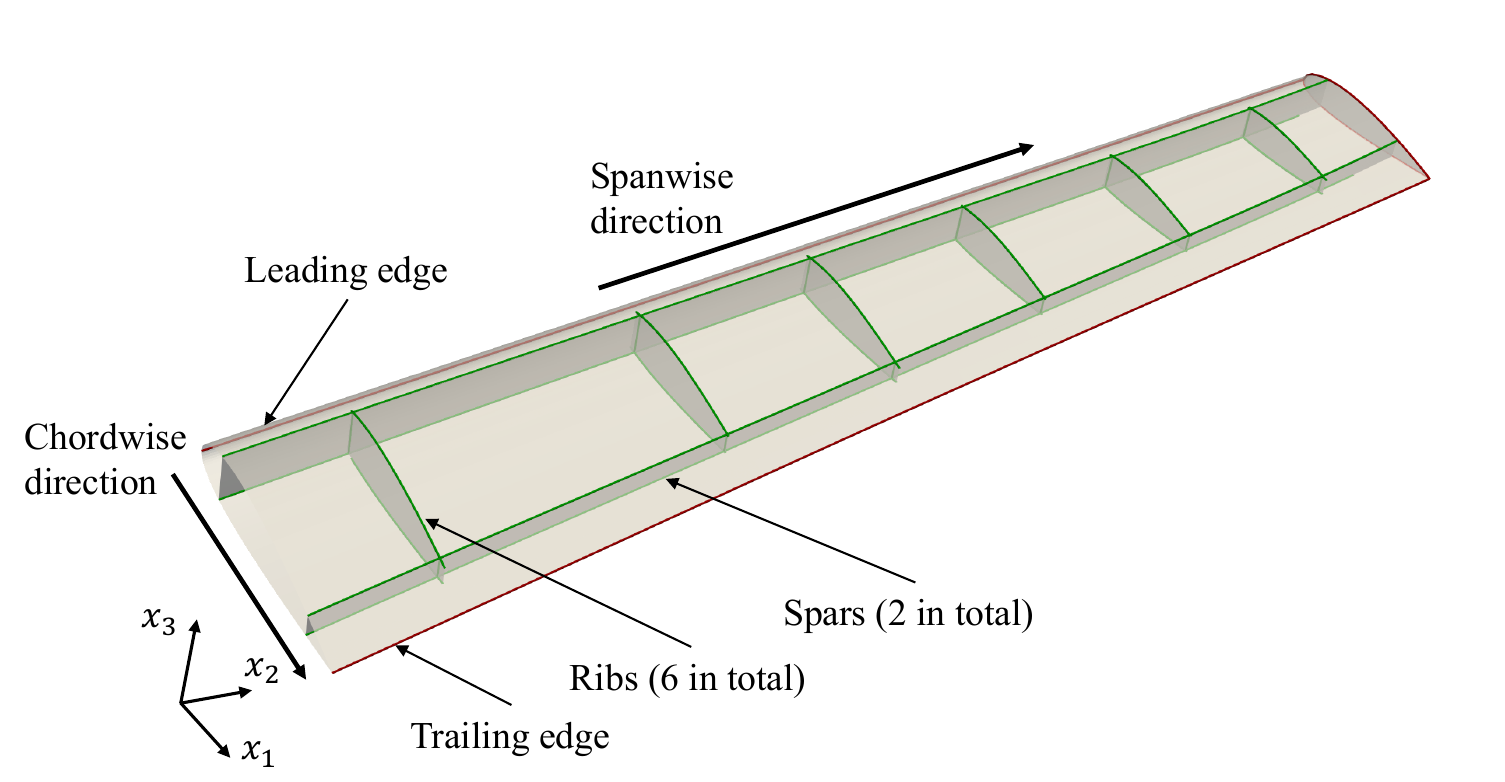}
    \caption{The CAD geometry of an eVTOL aircraft wing comprises 11 NURBS patches with 32 intersections. There are 28 movable intersections highlighted by green curves and 4 fixed intersections are indicated by red curves.}
    \label{fig:evtol-wing-geom}
\end{figure}

The proposed optimization scheme holds promise for enhancing the design of novel aerospace structures, where thin-walled structures are prevalent. We apply this shell shape optimization method with moving intersections to change the internal structures layout of an eVTOL aircraft wing, aiming to reduce the internal energy of the wing. The CAD geometry of the wing is depicted in Figure \ref{fig:evtol-wing-geom}, demonstrating the initial design created using the open-source software OpenVSP \cite{Gloudemans1996, Gloudemans2010, Fredericks2010, Hahn2010} developed by NASA. The wing geometry consists of 11 NURBS patches including 2 outer skins, 1 wing tip, 2 spars, and 6 ribs, where 32 intersections are detected in the wing geometry. Among the intersections, 4 of them are categorized as edge--edge intersections between outer surfaces or the wing tip, thus staying fixed throughout the optimization process and are marked with red curves in Figure \ref{fig:evtol-wing-geom}. The remaining intersections are either interior--interior, formed between ribs and spars, or interior--edge intersections, formed between outer surfaces and internal structures, and therefore can be moved during the optimization process. These movable intersections are distinguished by green curves.

\begin{figure}[!htb]\centering
    \includegraphics[width=0.9\textwidth]{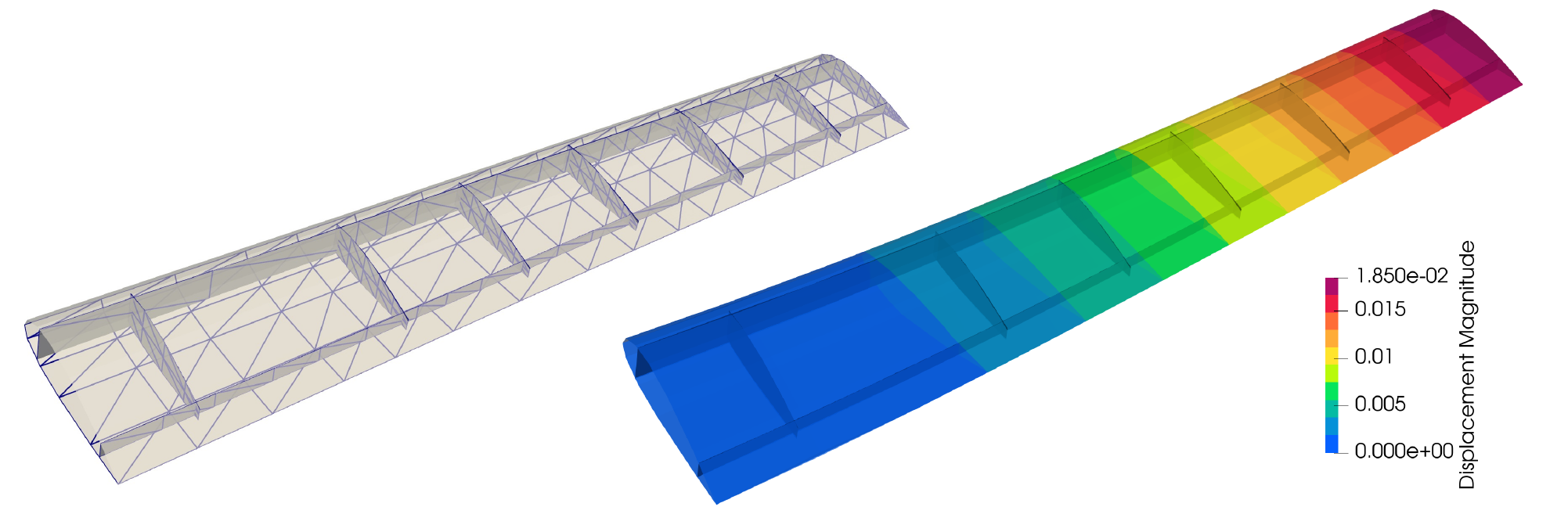}
    \caption{The initial eVTOL wing geometry discretization with NURBS basis functions, followed by the displacement result using the penalty-based non-matching coupling method for isogeometric Kirchhoff--Love shells. The displacement field has a unit of m and is scaled by a factor of 20 for visualization.}
    \label{fig:evtol-wing-opt-init}
\end{figure}

The baseline design of the wing geometry, isogeometrically discretized using cubic NURBS basis functions with a total of 2274 DoFs, is displayed in Figure  \ref{fig:evtol-wing-opt-init}, presenting the non-conforming discretizations between outer surfaces and internal structures. The wing geometry serves as the analysis model throughout the optimization process. A clamped boundary condition is applied at the wing root, while the lower outer surface is subjected to an upward distributed pressure of $500\ \rm{N}/\rm{m}^2$, simulating the cruise condition of the wing. For structural analysis, an isotropic elastic material model is employed for simplicity, with material properties corresponding to aluminum: Young's modulus $E=6.8\times 10^9\ \rm{Pa}$ and Poisson's ratio $\nu=0.35$. The wingspan is approximately 4.8 m and the chord is 1.2 m. All shell patches have a thickness of $3\ \rm{mm}$. The displacement field of the initial wing geometry is demonstrated in Figure \ref{fig:evtol-wing-opt-init}.

In this application problem, we first optimize the shape of the internal spars and ribs of the wing while keeping the outer surfaces unchanged to minimize the internal energy of the wing. For the first design scenario, we consider a rigid body approach, allowing only movement of the spars along the $x_1$ direction and the ribs along the $x_2$ direction. Thus, each spar and rib can only translate in a single direction, resulting in one associated design variable for the rigid body movement, totaling 8 design variables. However, due to the movement of the internal structures, the edges of these shell patches may deviate from the outer skins, requiring additional constraints to maintain the T-junctions. To achieve this, we utilize the constraints discussed in Section \ref{subsec:intersection-types}. Specifically, we employ a cubic NURBS curve with one knot span $[0,0,0,0,1,1,1,1]$ for each intersecting edge of the internal structures to preserve the T-junctions in the $x_3$ direction. Thus, each edge constraint requires 4 design variables. With 16 T-junctions in the wing geometry, this yields 64 design variables in the $x_3$ direction to the optimization problem, resulting in a total of 72 design variables for the rigid body shape optimization of the internal structures.

Furthermore, we introduce additional constraints to ensure that the spars remain within the envelope of the outer surfaces in the $x_1$ direction. This can be achieved by setting lower and upper bounds to the associated design variables. The front and rear edges of the ribs are maintained within 80\% and 15\% locations from the trailing edge to the leading edge in the chordwise direction. Meanwhile, a minimum distance of 0.4 m between the two spars at the wing root and a minimum distance of 0.5 m between adjacent ribs are imposed.

\begin{figure}[!htb]
    \centering
    \begin{subfigure}[t]{0.9\textwidth}
        \centering
        \includegraphics[width=\textwidth]{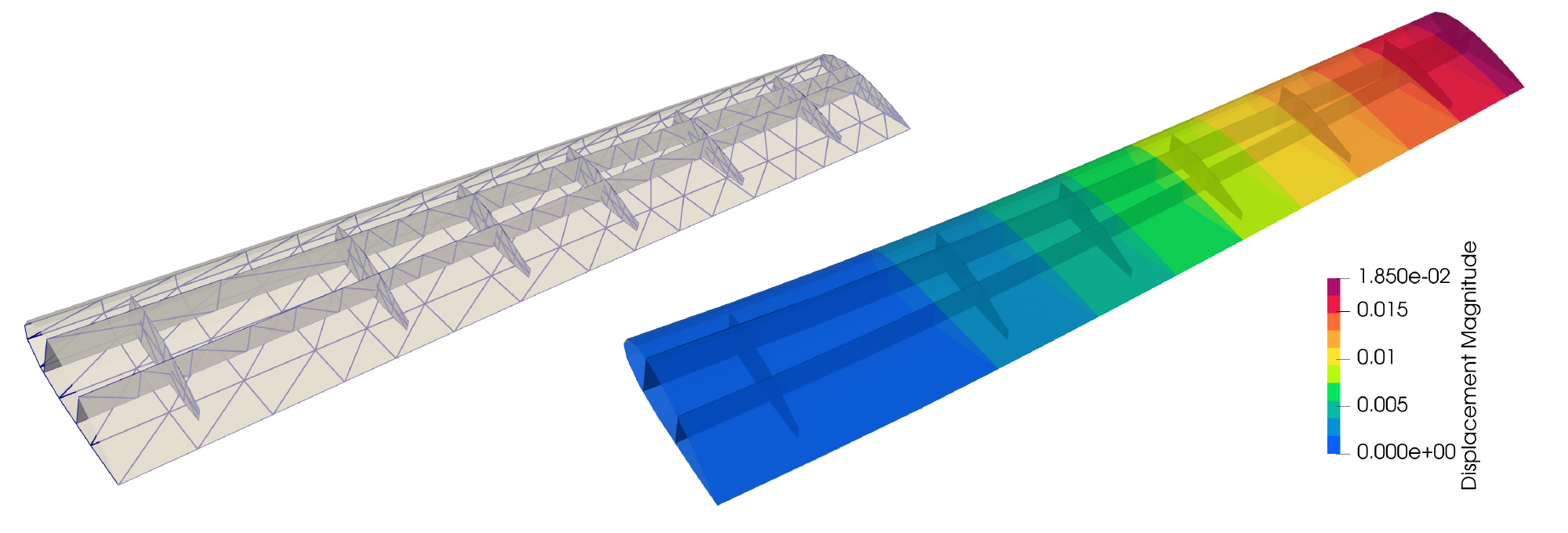}
        \caption{}
        \label{subfig:evtol-wing-opt-rs-rb}
    \end{subfigure}
    \hfill
    \begin{subfigure}[t]{0.9\textwidth}
        \centering
        \includegraphics[width=\textwidth]{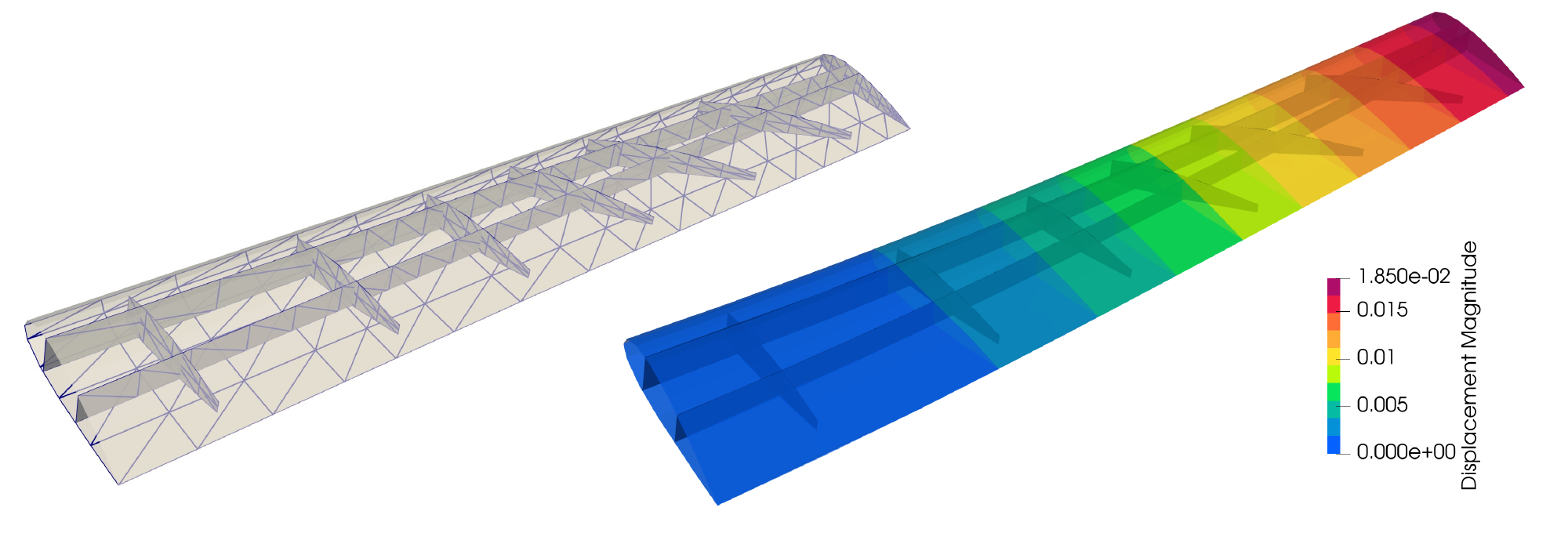}
        \caption{}
        \label{subfig:evtol-wing-opt-rs-sb}
    \end{subfigure}
    \caption{(a) The optimized geometry with rigid body translation for spars and ribs and associated displacements. (b) The optimized geometry with rigid body translation for spars and planar ribs and corresponding displacements. The displacement field is scaled by a factor of 20.}
    \label{fig:evtol-wing-opt-rs}
\end{figure}

With the problem setup described above, we proceed to perform rigid body shape optimization for the internal spars and ribs. In this design scenario, we utilize the SNOPT optimizer with tolerances of $10^{-3}$ for optimality and $10^{-4}$ for feasibility. Convergence is achieved after 22 iterations, and the optimal design and associated displacement magnitude are illustrated in Figure \ref{subfig:evtol-wing-opt-rs-rb}. In the optimized configuration, the two spars translate toward the center of the wing, reaching the minimum distance limit, while the ribs almost stay unchanged. The optimal design enhances the bending rigidity of the wing by moving the spars towards the center; a reasonable adjustment since the wing displacement is dominated by bending deformation given the geometry, loading and boundary conditions. On the other hand, the ribs, which are aligned parallel to the chordwise direction and are only allowed rigid body movement, exhibit negligible impact on reducing the bending deformation compared to the spars. The internal energy of the optimized wing geometry in Figure \ref{subfig:evtol-wing-opt-rs-rb} is 94.98\% of the baseline design in Figure \ref{fig:evtol-wing-opt-init}.

In the second design scenario, we enable the change of the 6 ribs for not only translation but also rotation in the $x_1$--$x_2$ plane, while maintaining the planar geometry. The front and rear edges of the ribs still stay along the lines of 80\% and 15\% of the distance from the trailing edge to the leading edge, respectively. Thus, the $x_1$ coordinates of the ribs are dependent on their $x_2$ coordinates. This adjustment increases the number of design variables in the $x_2$ direction to 2 for each rib, resulting in a total of 78 design variables. To prevent excessive rotation and elongation of the ribs, an additional set of constraints is introduced to ensure that the volume of each rib remains below 1.5 times the initial volume. The minimum distance between the two spars is set as 0.1 m. The same SNOPT optimizer and convergence criteria are employed for this optimization problem. The geometry and displacement associated with the optimized configuration after 46 iterations are demonstrated in Figure \ref{subfig:evtol-wing-opt-rs-sb}, where ribs increasingly tilt from the wing root to the wingtip, adding additional bending rigidity to the wing, particularly in regions where the two spars are in closer proximity. The internal energy of the optimized geometry is 94.84\% of the initial geometry, slightly lower than the first design scenario as we expected.

\begin{figure}[!htb]
    \centering
    \begin{subfigure}[t]{0.9\textwidth}
        \centering
        \includegraphics[width=\textwidth]{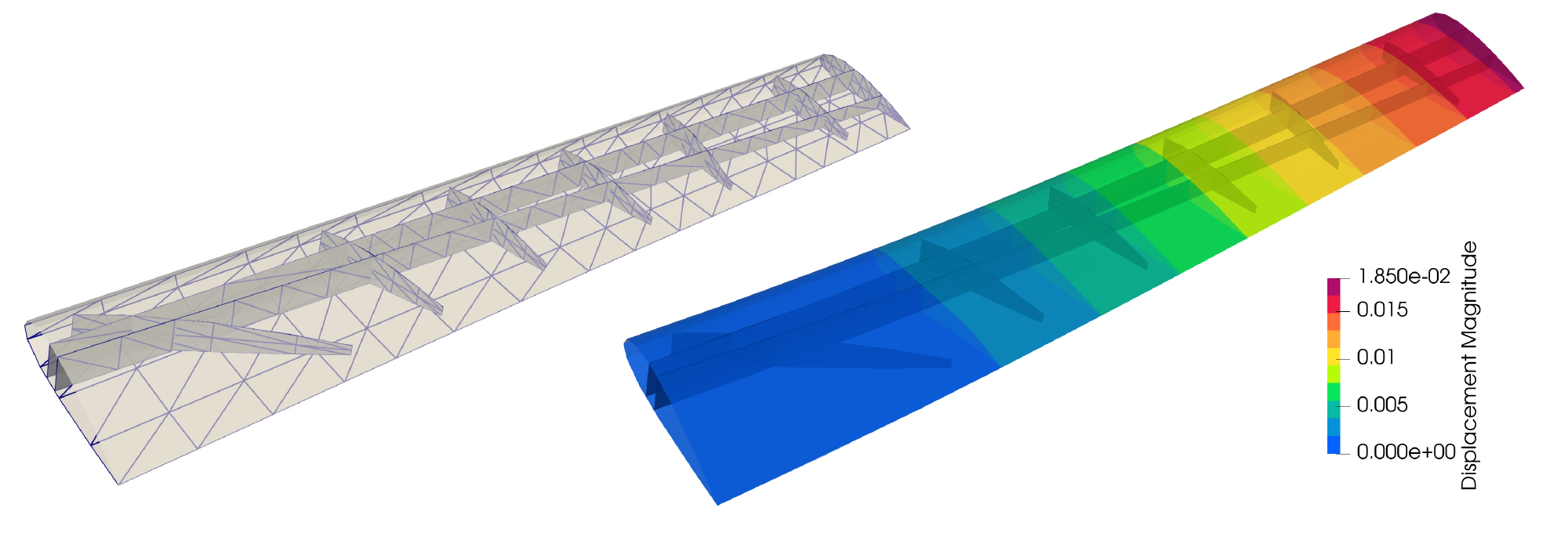}
        \caption{}
        \label{subfig:evtol-wing-opt-ss-sb}
    \end{subfigure}
    \hfill
    \begin{subfigure}[t]{0.9\textwidth}
        \centering
        \includegraphics[width=\textwidth]{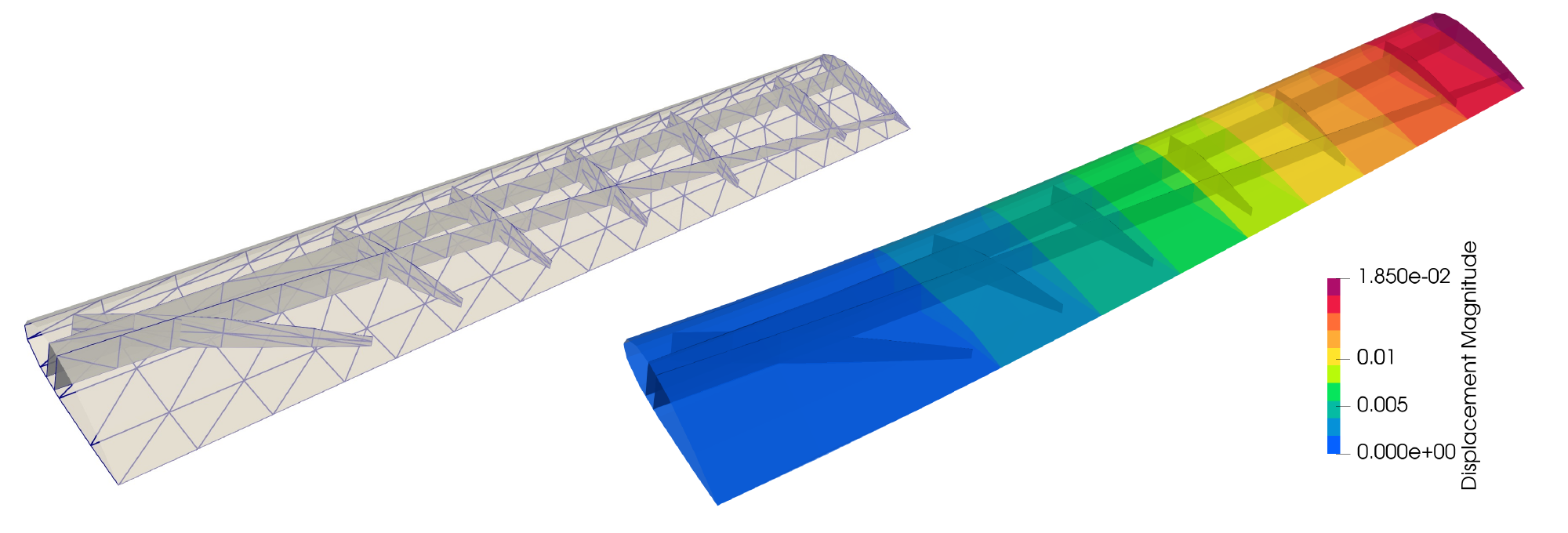}
        \caption{}
        \label{subfig:evtol-wing-opt-qs-sb}
    \end{subfigure}
    \caption{(a) The optimized geometry with planar spars and ribs and contour plot of the displacements. (b) The optimized geometry with quadratic spars and planar ribs and resulting displacements. The displacement field is scaled by a factor of 20.}
    \label{fig:evtol-wing-opt-ss}
\end{figure}

Considering that spars have more influence on the wing's bending deformation, we consider rotatable planar surfaces to model the spars in addition to the second design subjected to the same geometric and volume constraints. This enables all the internal structural members to translate and rotate in the third design scenario. Compared to the second case, one more design variable is introduced for each spar. Similar to the rib volume constraint, the volume of each spar is restricted to less than 1.5 times the initial volume. All other design conditions and optimization parameters remain unchanged. Convergence is achieved after 32 iterations, yielding an optimized geometry illustrated in Figure \ref{subfig:evtol-wing-opt-ss-sb}, where the corresponding contour plot of displacement is also included. Notably, the two spars move towards the center at the wing root and gradually split at the wingtip, while the ribs exhibit decreasing rotation from root to tip, in contrast to the second design case. Further, the results show a relation between the tilt degree of the ribs and the distance between the two spars; the closer the spars, the greater the rib tilt degrees. The internal energy of the third design is calculated to be 93.08\% of the initial geometry, representing a distinct improvement compared to the second design case in Figure \ref{subfig:evtol-wing-opt-rs-sb}.

To further improve the third design, we introduce an additional design variable to each spar in the $x_1$ direction so that the spars are described by quadratic B-splines. The B-splines feature a knot span of $[0,0,0,1,1,1]$ in the $x_1$--$x_2$ plane, allowing for out-of-plane curvature in addition to translation and rotation. This change increases the total number of design variables to 82. The optimizer converges in 25 iterations and yields an optimized geometry depicted in Figure \ref{subfig:evtol-wing-opt-qs-sb}, which follows a similar pattern to the third design case but incorporates curved spars. The internal energy of the optimized design is 93.04\% of the baseline design, slightly smaller than the third case. Considering other internal components in the wing and manufacturability, the third optimization scenario emerges as a more practical design, exhibiting a 6.92\% internal energy reduction compared to the baseline design. However, depending on the design conditions, the fourth design provides a nontraditional internal structure that may prove beneficial for other types of stiffened thin-wall structures.

The application to the eVTOL wing demonstrates the effectiveness of the proposed approach in shape optimization of complex shell structures involving large movement of surface intersections. Wing spars and ribs are reorganized based on design criteria to minimize the internal energy of the wing. 
As the dimension of the design space increases, the objective function converges to a smaller value as expected. In the optimized configuration, the edges of the internal structures align well with the outer surfaces, effectively retaining the T-junctions. Throughout the optimization process, all shell patches maintain analysis-suitable NRUBS surfaces without significant distortion in the discretization despite the large movement of intersections. The eVTOL wing structural components are entirely represented by the NURBS surfaces during the whole optimization process, ensuring precise shape updates in analysis and shape optimization under a streamlined workflow.

\begin{figure}[!htb]
    \centering
    \includegraphics[width=0.7\textwidth]{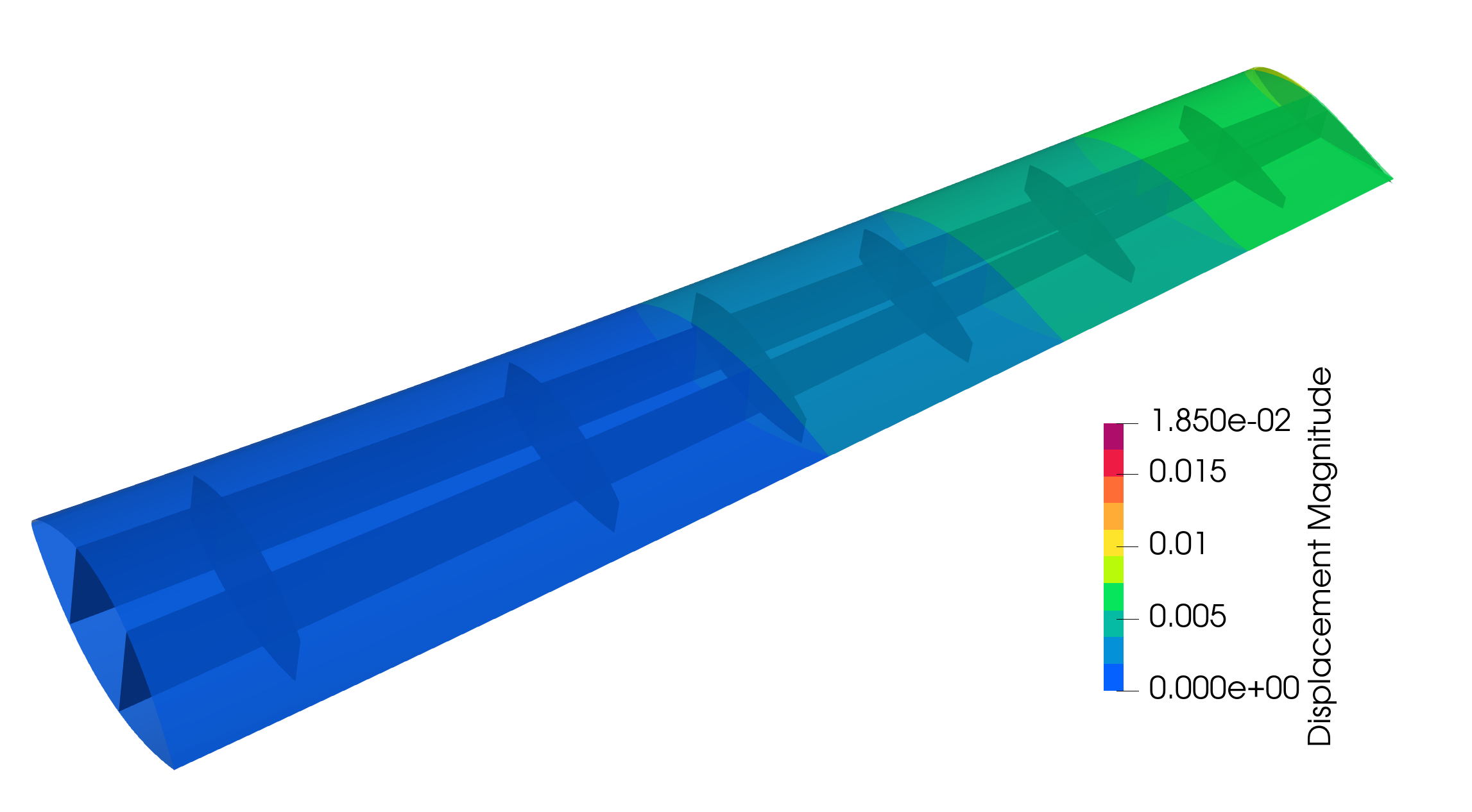}
    \caption{The displacement field of the optimized geometry incorporating updates of outer skins and rigid body translation of spars ribs and resulting displacements. The displacement field is scaled by a factor of 20.}
    \label{fig:evtol-wing-opt-outer}
\end{figure}

Finally, in addition to optimizing only the internal structures, considering shape optimization for both internal structures and outer skins simultaneously can achieve significant improvements for the eVTOL wing but yield a more challenging optimization problem. The internal structures must maintain the T-junctions along with the changing outer surfaces while searching for the optimal position. In this design scenario, we consider rigid body translation for the internal ribs and spars. The design models of the lower and upper outer skins are described using surfaces with a linear NURBS in the spanwise direction and a cubic NURBS in the chordwise direction, with knot vectors $[0,0,1,1]$ and $[0,0,0,0,1,1,1,1]$, respectively. To simplify the problem, we fix the shape of the edges of each outer skin except for the edge at the wing root and only change the vertical coordinates of the outer surfaces’ control points. This choice of design model introduces 8 new design variables for each outer skin surface. The upper bound of 2.9 m and lower bound of 3.3 m are set for the control points of the outer surfaces. The same boundary and loading conditions are applied to the eVTOL wing. The optimized geometry and associated displacement contour are depicted in Figure \ref{fig:evtol-wing-opt-outer}, where the internal energy of the wing is reduced by 64.54\%. The wing root expands in the vertical direction, and the spars move towards the center simultaneously to enhance the bending rigidity of the wing. Meanwhile, the longer edges of the ribs and spars still maintain the T-junctions with the lower and upper surfaces despite their shape updates. Due to the use of a simple distributed load, the control points of the outer surfaces reach the specified limits to maximize the support of the wing, as expected. Aerodynamic solvers and appropriate aero-structural coupling methods are required in future work to obtain a more realistic design.

\section{Conclusion}\label{sec:conclusion}
This paper introduces a novel shape optimization approach for shell structures composed of multiple NURBS patches. This method allows relative movement between shell patches without compromising mesh quality during shape updates, thereby enabling moving intersections. Shell patches are modeled using Kirchhoff--Love theory and discretized isogeometrically in the analysis, and intersecting shell patches are coupled using a penalty-based method. To achieve the moving intersections during design optimization, partial derivatives of the penalty energy with respect to surface intersections' parametric locations are formulated, along with an implicit relation between NURBS surfaces' control points and surface intersections' parametric locations. Standard benchmark problems are employed to validate this shape optimization approach, demonstrating that the optimized solutions closely match the reference solutions. Furthermore, we apply the proposed approach to adjust the layout of the internal spars and ribs of an eVTOL aircraft wing, resulting in nontraditional wing designs aimed at minimizing the internal energy of the wing.

In the optimization workflow, an initial CAD geometry of a shell structure is imported into the framework, with the control points of the CAD geometry treated as design variables. Their coordinates are updated directly using the optimization algorithm. The CAD geometry with updated control points is then used in the analysis framework directly to evaluate the structural response and total derivatives without additional steps. 

Integrating IGA into shell shape optimization presents notable advantages. The direct analysis based on CAD geometries in IGA naturally bridges the gap between the design model and analysis model within the optimization loop without geometric errors. The coarse design model is employed to update the shape of the geometry, while the fine analysis model is used for structural analysis. Shape modifications are directly applied to the coarse design model, and the structural response of the updated geometry is evaluated using the refined analysis model. As such, the dimension of the design space can be significantly reduced. This shape optimization workflow for non-matching shells is significantly simplified and accelerates the conceptual design of complex shell structures. The optimal designs provide valuable insights for the development of innovative shell structures.

\section*{Acknowledgements}\label{sec:Acknowledgements}
H. Zhao was supported by NASA grant number 80NSSC21M0070 while preparing the original submission. We thank Dr. David Kamensky at the University of California San Diego for helpful discussions on partial derivative computation algorithms.

\appendix
\section{Partial derivatives of the non-matching residual} \label{app:derivative-drpen-dxi}

The block matrices in \eqref{eq:expression-multi-patch-drcoudxi} are obtained through the application of the chain rule and taking derivatives for the NURBS basis functions. The formulation of the first diagonal block is detailed as follows,
\begin{equation}
    \begin{aligned}
    \partial_{\tilde{\bm{\xi}}^{\text{A}}} \mathbf{R}^{\text{A}}_{\text{pen}} 
    &= (\partial_{\tilde{\bm{\xi}}^{\text{A}}} \hat{\mathbf{R}}^{\text{A}}(\tilde{\bm{\xi}}^{\text{A}}))^{\mathrm{T}(2,3,1)} \, \tilde{\mathbf{R}}^{\text{A}}_{\text{pen}} 
    + (\hat{\mathbf{R}}^{\text{A}}(\tilde{\bm{\xi}}^{\text{A}}))^{\mathrm{T}} \, \partial_{\tilde{\bm{\xi}}^{\text{A}}} \tilde{\mathbf{R}}^{\text{A}}_{\text{pen}} \\[6pt]
    &+ (\partial_{\tilde{\bm{\xi}}^{\text{A}}} \hat{\mathbf{R}}^{\text{A}},_{\bm{\xi}}(\tilde{\bm{\xi}}^{\text{A}}))^{\mathrm{T}(2,3,1)} \, \tilde{\mathbf{R}}^{\text{A}\bm{\xi}}_{\text{pen}} 
    + (\hat{\mathbf{R}}^{\text{A}},_{\bm{\xi}}(\tilde{\bm{\xi}}^{\text{A}}))^{\mathrm{T}} \, \partial_{\tilde{{\bm{\xi}}}^{\text{A}}} \tilde{\mathbf{R}}^{\text{A}\bm{\xi}}_{\text{pen}} \text{ ,} \label{appeq:partial-rpena-partial-xia}
    \end{aligned}
\end{equation}
where $\tilde{\mathbf{R}}^{\text{A}}_{\text{pen}}=\partial_{\tilde{\mathbf{d}}^{\text{A}}} W^{\text{AB}}_{\text{pen}}$ and $\tilde{\mathbf{R}}^{\text{A}\bm{\xi}}_{\text{pen}}=\partial_{\tilde{\mathbf{d}}_{\bm{\xi}}^{\text{A}}} W^{\text{AB}}_{\text{pen}}$ are residual vectors of the penalty energy on the quadrature mesh $\tilde{\Omega}$ at location $\tilde{\bm{\xi}}^{\text{A}}$. $\partial_{\tilde{\bm{\xi}}^{\text{A}}} \hat{\mathbf{R}}^{\text{A}}(\tilde{\bm{\xi}}^{\text{A}})$ and $\partial_{\tilde{\bm{\xi}}^{\text{A}}} \hat{\mathbf{R}}^{\text{A}},_{\bm{\xi}}(\tilde{\bm{\xi}}^{\text{A}})$ are 3D arrays involve the first order and second order derivatives of the NURBS basis functions, where the transpose $\mathrm{T}(2,3,1)$ indicate switching the axes of the 3D array from $(1,2,3)$ to $(2,3,1)$. 
Partial derivative $\partial_{\tilde{\bm{\xi}}^{\text{A}}} \tilde{\mathbf{R}}^{\text{A}}_{\text{pen}}$ in \eqref{appeq:partial-rpena-partial-xia} is formulated as
\begin{equation}
    \begin{aligned}
        \partial_{\tilde{\bm{\xi}}^{\text{A}}} \tilde{\mathbf{R}}^{\text{A}}_{\text{pen}} 
        &= \partial_{\tilde{\mathbf{d}}^{\text{A}}} \tilde{\mathbf{R}}^{\text{A}}_{\text{pen}}\, \partial_{\tilde{\bm{\xi}}^{\text{A}}} \tilde{\mathbf{d}}^{\text{A}} 
        + \partial_{\tilde{\mathbf{d}}_{\bm{\xi}}^{\text{A}}} \tilde{\mathbf{R}}^{\text{A}}_{\text{pen}} \, \partial_{\tilde{\bm{\xi}}^{\text{A}}} \tilde{\mathbf{d}}_{\bm{\xi}}^{\text{A}} \\[6pt]
        &+ \partial_{\tilde{\mathbf{P}}^{\text{A}}} \tilde{\mathbf{R}}^{\text{A}}_{\text{pen}}\, \partial_{\tilde{\bm{\xi}}^{\text{A}}} \tilde{\mathbf{P}}^{\text{A}} 
        + \partial_{\tilde{\mathbf{P}}_{\bm{\xi}}^{\text{A}}} \tilde{\mathbf{R}}^{\text{A}}_{\text{pen}} \, \partial_{\tilde{\bm{\xi}}^{\text{A}}} \tilde{\mathbf{P}}_{\bm{\xi}}^{\text{A}} \text{ ,} \label{appeq:partial-rpena-tilda-partial-xia}
    \end{aligned}
\end{equation}
where the first component in each term, e.g., $\partial_{\tilde{\mathbf{d}}^{\text{A}}} \tilde{\mathbf{R}}^{\text{A}}_{\text{pen}}$, can be derived from the penalty energy \eqref{eq:penatly-virtual-work}. For the second components, such as $\partial_{\tilde{\bm{\xi}}^{\text{A}}} \tilde{\mathbf{d}}^{\text{A}}$, it is computed via the interpolation matrix,
\begin{equation}
    \begin{aligned}
        \partial_{\tilde{\bm{\xi}}^{\text{A}}} \tilde{\mathbf{d}}^{\text{A}} = \partial_{\tilde{\bm{\xi}}^{\text{A}}} (\hat{\mathbf{R}}^{\text{A}}(\tilde{\bm{\xi}}^{\text{A}})\, \mathbf{d}^{\text{A}}) = (\partial_{\tilde{\bm{\xi}}^{\text{A}}} \hat{\mathbf{R}}^{\text{A}}(\tilde{\bm{\xi}}^{\text{A}}))^{\mathrm{T}(1,3,2)}\,  \mathbf{d}^{\text{A}} \text{ .} \label{appeq:partial-da-partial-xia}
    \end{aligned}
\end{equation}
Likewise, the other component, $\partial_{\tilde{\bm{\xi}}^{\text{A}}} \tilde{\mathbf{d}}_{\bm{\xi}}^{\text{A}}$, involves the first derivative of the shell displacement is calculated as
\begin{equation}
    \begin{aligned}
        \partial_{\tilde{\bm{\xi}}^{\text{A}}} \tilde{\mathbf{d}}_{\bm{\xi}}^{\text{A}} = (\partial_{\tilde{\bm{\xi}}^{\text{A}}} \hat{\mathbf{R}}^{\text{A}},_{\bm{\xi}}(\tilde{\bm{\xi}}^{\text{A}}))^{\mathrm{T}(1,3,2)}\,  \mathbf{d}^{\text{A}} \text{ .} \label{appeq:partial-daxi-partial-xia}
    \end{aligned}
\end{equation}
Replacing $\mathbf{d}^{\text{A}}$ with $\mathbf{P}^{\text{A}}$ in \eqref{appeq:partial-da-partial-xia} and \eqref{appeq:partial-daxi-partial-xia}, we can obtain $\partial_{\tilde{\bm{\xi}}^{\text{A}}} \tilde{\mathbf{P}}^{\text{A}}$ and $\partial_{\tilde{\bm{\xi}}^{\text{A}}} \tilde{\mathbf{P}}_{\bm{\xi}}^{\text{A}}$, and therefore, $\partial_{\tilde{\bm{\xi}}^{\text{A}}} \tilde{\mathbf{R}}^{\text{A}}_{\text{pen}}$. The same derivation can be applied to obtain $\partial_{\tilde{{\bm{\xi}}}^{\text{A}}} \tilde{\mathbf{R}}^{\text{A}\bm{\xi}}_{\text{pen}}$. Substituting $\partial_{\tilde{\bm{\xi}}^{\text{A}}} \mathbf{R}^{\text{A}}_{\text{pen}}$ and $\partial_{\tilde{{\bm{\xi}}}^{\text{A}}} \tilde{\mathbf{R}}^{\text{A}\bm{\xi}}_{\text{pen}}$ into \eqref{appeq:partial-rpena-partial-xia}, we can get the first diagonal block matrix of the partial derivative $\partial_{\tilde{\bm{\xi}}} \mathbf{R}_{}$.

The off-diagonal block has a similar formulation to the diagonal block,
\begin{align}
    \partial_{\tilde{\bm{\xi}}^{\text{B}}} \mathbf{R}^{\text{A}}_{\text{pen}} 
    = (\hat{\mathbf{R}}^{\text{A}}(\tilde{\bm{\xi}}^{\text{A}}))^{\mathrm{T}} \, \partial_{\tilde{\bm{\xi}}^{\text{B}}} \tilde{\mathbf{R}}^{\text{A}}_{\text{pen}}
    + (\hat{\mathbf{R}}^{\text{A}},_{\bm{\xi}}(\tilde{\bm{\xi}}^{\text{A}}))^{\mathrm{T}} \, \partial_{\tilde{{\bm{\xi}}}^{\text{B}}} \tilde{\mathbf{R}}^{\text{A}\bm{\xi}}_{\text{pen}} \text{ .} \label{appeq:partial-rpena-partial-xib}
\end{align}
With the diagonal and off-diagonal blocks, we can obtain the partial derivative in \eqref{eq:expression-multi-patch-drcoudxi}. 

In \eqref{appeq:partial-rpena-partial-xia} and \eqref{appeq:partial-daxi-partial-xia}, it becomes evident that the second-order derivative of the NURBS basis functions for the shell patch is included. Consequently, the $C^1$ continuity is required if the evaluation point is located at the element boundary. This requirement is inherently fulfilled by the NURBS basis functions, thereby highlighting the advantages of IGA in the application of design optimization.

\section{Partial derivatives of the implicit intersection representation} \label{app:derivative-drl-dxi}

To compute the analytical total derivative of the non-matching shell optimization problem, partial derivatives of the implicit intersection representation \eqref{eq:expression-multi-patch-p-xi-residual} with respect to $\mathbf{P}$ and $\tilde{\bm{\xi}}$ are required. The partial derivative of the implicit residual vector with respect to intersections' parametric coordinates, $\partial_{\tilde{\bm{\xi}}} \mathbf{R}_{\mathcal{L}}$, states as
\begin{align}
    \partial_{\tilde{\bm{\xi}}} \mathbf{R}_{\mathcal{L}} = \begin{bmatrix}
        \mathbf{E}^{\text{A}} & \quad \mathbf{E}^{\text{B}}\\
        \mathbf{F}^{\text{A}} & \quad \mathbf{0} \\ 
        \multicolumn{2}{c}{\mathbf{v}_{r}} \\ 
        \multicolumn{2}{c}{\mathbf{v}_{s}}
    \end{bmatrix} \text{ ,} \label{appeq:partial-rl-partial-xi}
\end{align}
where 
\begin{equation}
\begin{aligned}
    \mathbf{E}_{ik}^{\text{A}} &= \begin{cases}
        \begin{bmatrix}
            \hat{\mathbf{R}}^{\text{A}},_{\xi_1}(\tilde{\bm{\xi}}_i^{\text{A}}) \mathbf{P}^{\text{A}} & \hat{\mathbf{R}}^{\text{A}},_{\xi_2}(\tilde{\bm{\xi}}_i^{\text{A}}) \mathbf{P}^{\text{A}}
        \end{bmatrix}\text{ ,} \quad &\text{ if } i=k\\
        \mathbf{0} \text{ ,} \quad &\text{ otherwise}
    \end{cases} \quad \text{ and } \\[6pt]
    \mathbf{E}_{ik}^{\text{B}} &= \begin{cases}
        -\begin{bmatrix}
            \hat{\mathbf{R}}^{\text{B}},_{\xi_1}(\tilde{\bm{\xi}}_i^{\text{B}}) \mathbf{P}^{\text{B}} & \hat{\mathbf{R}}^{\text{B}},_{\xi_2}(\tilde{\bm{\xi}}_i^{\text{B}}) \mathbf{P}^{\text{B}}
        \end{bmatrix}\text{ ,} \quad &\text{if } i=k\\
        \mathbf{0}\text{ ,} \quad &\text{ otherwise}
    \end{cases}\text{ ,}
\end{aligned}
\end{equation}
for $i=\{1,2,\ldots,m\}$ and $k=\{1,2,\ldots,m\}$. Matrices $\mathbf{E}^{\text{A}}$ and $\mathbf{E}^{\text{B}}$ have sizes of $(m\times sd)\times (m\times pd)$. Each entry, $\mathbf{E}_{ik}^{\text{A}}$ or $\mathbf{E}_{ik}^{\text{B}}$, is a block matrix with a size of $sd\times pd$. And 
\begin{equation}
\begin{aligned}
    \mathbf{F}_{jk}^{\text{A}} &= \begin{cases}
        2\left(\hat{\mathbf{R}}^{\text{A}}(\tilde{\bm{\xi}}_j^{\text{A}}) \mathbf{P}^{\text{A}} - \hat{\mathbf{R}}^{\text{A}}(\tilde{\bm{\xi}}_{j-1}^{\text{A}}) \mathbf{P}^{\text{A}} \right)^{\mathrm{T}} \begin{bmatrix}
            \hat{\mathbf{R}}^{\text{A}},_{\xi_1}(\tilde{\bm{\xi}}_{j-1}^{\text{A}}) \mathbf{P}^{\text{A}} & \hat{\mathbf{R}}^{\text{A}},_{\xi_2}(\tilde{\bm{\xi}}_{j-1}^{\text{A}}) \mathbf{P}^{\text{A}}
        \end{bmatrix}\text{ ,} \quad &\text{ if } k=j-1\\[8pt]
        -2\left(\hat{\mathbf{R}}^{\text{A}}(\tilde{\bm{\xi}}_{j+1}^{\text{A}}) \mathbf{P}^{\text{A}} - \hat{\mathbf{R}}^{\text{A}}(\tilde{\bm{\xi}}_{j-1}^{\text{A}}) \mathbf{P}^{\text{A}} \right)^{\mathrm{T}} \begin{bmatrix}
            \hat{\mathbf{R}}^{\text{A}},_{\xi_1}(\tilde{\bm{\xi}}_j^{\text{A}}) \mathbf{P}^{\text{A}} & \hat{\mathbf{R}}^{\text{A}},_{\xi_2}(\tilde{\bm{\xi}}_j^{\text{A}}) \mathbf{P}^{\text{A}}
        \end{bmatrix}\text{ ,} \quad &\text{ if } k=j\\[8pt]
        2\left(\hat{\mathbf{R}}^{\text{A}}(\tilde{\bm{\xi}}_{j+1}^{\text{A}}) \mathbf{P}^{\text{A}} - \hat{\mathbf{R}}^{\text{A}}(\tilde{\bm{\xi}}_{j}^{\text{A}}) \mathbf{P}^{\text{A}} \right)^{\mathrm{T}} \begin{bmatrix}
            \hat{\mathbf{R}}^{\text{A}},_{\xi_1}(\tilde{\bm{\xi}}_{j+1}^{\text{A}}) \mathbf{P}^{\text{A}} & \hat{\mathbf{R}}^{\text{A}},_{\xi_2}(\tilde{\bm{\xi}}_{j+1}^{\text{A}}) \mathbf{P}^{\text{A}}
        \end{bmatrix}\text{ ,} \quad &\text{ if } k=j+1\\[8pt]
        \mathbf{0} \text{ ,} \quad &\text{ otherwise}
    \end{cases} \text{ ,}
\end{aligned}
\end{equation}
for $j\in \{2,3,\ldots,m-1\}$ and $k\in \{1,2,\ldots,m\}$. $\mathbf{F}^{\text{A}}$ has dimensions of $(m-2)\times (m\times pd)$, with each entry $\mathbf{F}_{jk}^{\text{A}}$ being of size of $1\times pd$. Additionally, $\mathbf{v}_{r}$ is a row vector consisting of zero values except for the $r$-th entry, which is set to 1. And $\mathbf{v}_{s}$ possesses these same properties.

The partial derivative of the residual vector with respect to shell patches' control points $\partial_{\mathbf{P}} \mathbf{R}_{\mathcal{L}}$ is expressed as
\begin{align}
    \partial_{\mathbf{P}} \mathbf{R}_{\mathcal{L}} = \begin{bmatrix}
        \hat{\mathbf{R}}^{\text{A}}(\tilde{\bm{\xi}}_i^{\text{A}}) & -\hat{\mathbf{R}}^{\text{B}}(\tilde{\bm{\xi}}_i^{\text{B}}) \\ 
        \mathbf{G}^{\text{A}}_j & \mathbf{0}\\ 
        \multicolumn{2}{c}{\mathbf{0}} \\ 
        \multicolumn{2}{c}{\mathbf{0}}
    \end{bmatrix} \text{ ,} \label{appeq:partial-rl-partial-p}
\end{align}
for $i\in \{1,2,\ldots,m\}$ and $j\in \{2,3,\ldots,m-1\}$. And $\mathbf{G}^{\text{A}}_j$ is a row vector and has the following definition
\begin{equation}
    \begin{aligned}
        \mathbf{G}^{\text{A}}_j = 2{\mathbf{P}^{\text{A}}}^{\mathrm{T}}&\left[\left(\hat{\mathbf{R}}^{\text{A}}(\tilde{\bm{\xi}}_{j+1}^{\text{A}}) - \hat{\mathbf{R}}^{\text{A}}(\tilde{\bm{\xi}}_{j}^{\text{A}})  \right)^{\mathrm{T}} \left( \hat{\mathbf{R}}^{\text{A}}(\tilde{\bm{\xi}}_{j+1}^{\text{A}}) - \hat{\mathbf{R}}^{\text{A}}(\tilde{\bm{\xi}}_{j}^{\text{A}}) \right) \right. \\[6pt]
        &\left.- \left(\hat{\mathbf{R}}^{\text{A}}(\tilde{\bm{\xi}}_j^{\text{A}}) - \hat{\mathbf{R}}^{\text{A}}(\tilde{\bm{\xi}}_{j-1}^{\text{A}}) \right)^{\mathrm{T}} \left( \hat{\mathbf{R}}^{\text{A}}(\tilde{\bm{\xi}}_j^{\text{A}}) - \hat{\mathbf{R}}^{\text{A}}(\tilde{\bm{\xi}}_{j-1}^{\text{A}}) \right)\right] \text{ .}
    \end{aligned}
\end{equation}

Consequently, we can obtain the total derivative of the intersections' parametric coordinates with respect to shell patches' control points by substituting \eqref{appeq:partial-rl-partial-xi} and \eqref{appeq:partial-rl-partial-p} into \eqref{eq:expression-multi-patch-dxidp}.

\bibliographystyle{unsrt}
\bibliography{main}

\end{document}